\documentclass[letterpaper, 10 pt, conference]{ieeeconf}

\IEEEoverridecommandlockouts
\usepackage[LGR,T1]{fontenc}
\usepackage[latin9]{inputenc}
\usepackage{float}
\usepackage{amsthm}
\usepackage{amsbsy}
\usepackage{amstext}
\usepackage{amssymb}
\usepackage{amsmath}
\usepackage{esint}
\usepackage{graphicx}

\usepackage{color}

\DeclareFontEncoding{LGR}{}{}
\DeclareTextSymbol{\~}{LGR}{126}

\floatstyle{ruled}
\newfloat{algorithm}{tbp}{loa}
\providecommand{\algorithmname}{Algorithm}
\floatname{algorithm}{\protect\algorithmname}

\newtheorem{thm}{\protect\theoremname}
\newtheorem{rem}[thm]{\protect\remarkname}
\newtheorem{defn}[thm]{\protect\definitionname}

\providecommand{\definitionname}{Definition}
\providecommand{\remarkname}{Remark}
\providecommand{\theoremname}{Theorem}

\begin{document}

\title{Linear Hamilton Jacobi Bellman Equations in High Dimensions}

\author{Matanya B. Horowitz, Anil Damle, Joel W. Burdick
\thanks{Anil Damle is supported by NSF Fellowship DGE-1147470.
        The corresponding author is available at {\tt\small mhorowit@caltech.edu}
}
}

\maketitle
\thispagestyle{empty}
\pagestyle{empty}

\begin{abstract}
The Hamilton Jacobi Bellman Equation (HJB) provides the globally optimal solution to large classes of control problems. Unfortunately, this generality comes at a price, the calculation of such solutions is typically intractible for systems with more than moderate state space size due to the \emph{curse of dimensionality}. This work combines recent results in the structure of the HJB, and its reduction to a linear Partial Differential Equation (PDE), with methods based on low rank tensor representations, known as a separated representations, to address the curse of dimensionality. The result is an algorithm to solve optimal control problems which scales linearly with the number of states in a system, and is applicable to systems that are nonlinear with stochastic forcing in finite-horizon, average cost, and first-exit settings. The method is demonstrated on inverted pendulum, VTOL aircraft, and quadcopter models, with system dimension two, six, and twelve respectively.
\end{abstract}

\section{Introduction}
\label{sec:Introduction}
The Hamilton Jacobi Bellman (HJB) equation is central to control theory, yielding the optimal solution to general problems specified by known dynamics and a specified cost functional. Given the assumption of quadratic cost on the control input, it is well known that the HJB reduces to a particular Partial Differential Equation (PDE) \cite{Fleming:2006tl}. While powerful, this reduction is not commonly used as the PDE is of second order,  is nonlinear, and examples exist where the problem may not have a solution in a classical sense \cite{Crandall:1992ta}. Furthermore,
each state of the system appears as another dimension of the PDE, giving rise to the \emph{curse of dimensionality} \cite{Bellman:1962vr}. Since the number of degrees of freedom required to solve the optimal control problem grows exponentially with dimension, the problem becomes intractable for systems with all but modest dimension. 

In the last decade researchers have found that under certain, fairly
non-restrictive, structural assumptions, the HJB may be transformed
into a linear PDE, see, \emph{e.g.}, \cite{Todorov:2009je} and \cite{Kappen:2005bn}, with an interesting
analogue in the discretized domain of Markov Decision Processes (MDP)
\cite{Todorov:2006tq}. The implications for this discovery are numerous,
and research has only begun to tap into the computational benefits
\cite{dvijotham2012linearly} and \cite{Horowitz:2014tu}.
The work presented here is a continuation of this theme, and uses
the linearity of this particular form of the HJB PDE to push the computational
boundaries of the HJB.

Our method relies on recent work in Separated Representations (SR) \cite{Beylkin:2005kh}, which have recently emerged as a method to solve a number of problems in machine learning and the numerical solution of PDEs with complexity that scales \emph{linearly} with dimension, bypassing the curse of dimensionality. The central idea of this paper is to approximate the solution, and its associated operators, by a low rank tensor. If the problem's components can be adequately modeled in this regime, then the complexity grows with the rank of the approximation, rather than the dimensionality. For many problems of interest this proves to be a valid modeling assumption.

\section{Related Work}
\label{sec:Related_work}
We combine two previously disjoint threads of research. The first is the study
of HJB equations, while the second is the study of high dimensional tensors and their approximations.

\subsection{Linearly Solvable Stochastic Optimal Control}
The study of linearly solvable Stochastic Optimal Control (SOC) problems has developed along two lines of investigation. One is that of Linear MDPs \cite{Todorov:2009je}, in which a control design problem may be solved via a linear set of equations given several structural assumptions. By taking the continuous limit of the discretization, a linear PDE is obtained. In another line of work begun by Kappen \cite{Kappen:2005bn} the same linear PDE has been found through a particular transformation of the HJB. While the linearity of the HJB
provides computational benefits in terms of the numerical techniques
available, the curse of dimensionality prevents existing techniques
from scaling to realistic problems of interest. This has been addressed
through the use of the Path Integral techniques \cite{Theodorou:2010vd,Theodorou:2011uz}, which rely on Monte Carlo techniques via the Feynman-Kac Lemma. The solution of the Linear PDE at an individual point in the state space is solved via an ensemble of Brownian motions. While
these sampling techniques are formally independent of state space
dimension, they may nonetheless be computationally expensive and typically
only provide guarantees in the asymptotic regime. Furthermore, as the solution is only valid at an individual point, these solutions must be solved for all anticipated states that will occur over the course of a trajectory. These techniques have nonetheless shown notable success \cite{Stulp:2012jb}.

\subsection{Nonlinear Hamilton Jacobi Bellman Equations in Control}
HJB equations have arisen as objects of study in the research of Control
Lyapunov Functions. These functions generalize the notion of Lyapunov
stability, allowing for the control signal to be incorporated in the analysis of a system's stability \cite{Sontag:1983cq}. These techniques may be seen as a relaxation of the conditions that lead to optimality for the HJB, and thus existing CLF-synthesis techniques are in general suboptimal. This may be ameliorated by combining CLFs with Euler-Lagrange equations which arise from the Pontraygin maximum principle, typically recognized as Model Predictive Control or Receding Horizon Control \cite{Primbs:1999tt}. Indeed, many planning techniques may be understood as heuristic or approximate methods to solve the HJB \cite{lavalle2006planning, Horowitz:nav}.

Ours is not the first study of numerical techniques to solve the HJB
PDE directly. In \cite{Mitchell:2003vx} level-set
algorithms are used to solve variants of the HJB that relate to the
calculation of reachable sets. In \cite{Aguilar:2013er}, the authors uses high order
Taylor expansions to approximate the HJB directly. McEaney \cite{mceneaney2007curse} has developed another curse-of-dimensionality free method which relies on a max-plus expansion of the solution, with complexity that scales with the number of basis solutions, each requiring the solution to a Riccati system. 

Lasserre provides an alternative framework to solve the HJB by constraining the moments of the solution, producing bounds on the moments that can be made to converge through a hierarchy of optimization problems \cite{Lasserre:2008ej,Lasserre:2001fw}. These techniques seek to reduce the fundamental infinite dimensional linear program to a semidefinite program with a finite number of degrees of freedom by truncating the solution moments. This is currently an active and growing area of research \cite{majumdar2013control}. A similar construction of Lasserre hierarchies
was proposed in \cite{Horowitz:2014tu} using Sum of Squares techniques,
producing upper and lower bounds to the linear HJB. However, all fall
prey to the curse of dimensionality, save the work of McEaney, and in practice are limited to dimensions of five or lower. We argue that while the techniques presented herein do not yet apply to all optimal control problems at the moment, they apply to a broad class of systems, with astounding computational gains. We discuss how these techniques may in fact be applied without these assumptions in our concluding discussion.

Researchers have previously attacked the intractability of the HJB
through discretization of the system state space, creating an MDP. The curse
of dimensionality is mitigated in this context by parameterizing the
value function with a sparse set of basis, giving rise to Approximate
Dynamic Programming, or when the basis may change online Adaptive
Dynamic Programming (ADP) \cite{Bertsekas:2012uq}. These techniques have constraint sets that formally grow exponentially with dimensionality \cite{deFarias:2003wq}. Nonetheless, these techniques are the most popular method to deal with the curse of dimensionality, and have even been used to surpass human capabilities on complex time dependent games via synthesis with modern machine learning techniques \cite{Mnih:2013wp}. These methods are closest to ours in spirit, and our method could be seen as generating a sparse basis, as is desired in ADP, albeit ours is performed without recourse to an MDP, with the attendant constraints.

\subsection{High Dimensional Tensors}

Tensor approximations have historically been developed with the goal
of approximating high dimensional data, yielding rise to the framework
used here under the names CANDECOMP/PARAFAC \cite{Harshman:1970vk,Carroll:1970bw}. However,  \cite{Beylkin:2005kh} demonstrated that these approximation techniques were applicable to the linear systems describing discretized PDEs. This technique has been applied in several domains, including computational chemistry and quantum physics, among others \cite{Khoromskij:2012hq}. In particular, \cite{Sun2014} examines their use in the context of stationary Fokker-Planck equations. There are interesting connections between the fundamental task of these techniques,
approximating a tensor with one of lower rank, and convex relaxation
based methods \cite{Chandrasekaran:2012hl,Gandy:2011bh}. Unfortunately,
low rank tensor approximation is NP-hard in general, and an optimal
solution is not to be expected \cite{Hastad:1990tx}.

\section{The Linear Hamilton Jacobi Bellman Equation}
\label{sec:The_Linear_HJB}

We begin by constructing the value function, which captures the \emph{cost-to-go} from a given state. If such a quantity is known, the optimal action
may be chosen as that which follows the gradient of the value, bringing
the agent into the states with highest value over the remaining time
horizon. We define $x_{t}\in\mathbb{R}^{n}$ as the system state at time $t$,
control input $u_{t}\in\mathbb{R}^{m}$ , and dynamics that evolve
according to the equation
\begin{equation}
dx_{t}=\left(f\left(x_{t}\right)+G\left(x_{t}\right)u_{t}\right)dt+B\left(x_{t}\right)d\omega_{t}\label{eq:dynamics}
\end{equation}
on a compact domain $\Omega$, and where the expressions $f(x)$, $G(x)$,
$B(x)$ are assumed to be smoothly differentiable, but possibly nonlinear,
functions, and $\omega_{t}$ is Gaussian noise with covariance $\Sigma_\epsilon$.
The system has cost $r_{t}$ accrued at time $t$ according to
\begin{equation}
r\left(x_{t},u_{t}\right)=q\left(x_{t}\right)+\frac{1}{2}u_{t}^{T}Ru_{t}\label{eq:cost}
\end{equation}
where $q(x)$ is a state dependent cost. We require $q(x)\ge0$ for all $x$ in the problem domain and $R$ positive definite. The goal is to minimize the expectation of the cost functional 
\begin{equation}
J(x,u)=\phi_{T}\left(x_{T}\right)+\int_{0}^{T}r\left(x_{t},u_{t}\right)dt\label{eq:trajectory-cost}
\end{equation}
where $\phi_{T}$ represents a state-dependent terminal cost. The
solution to this minimization is known as the \emph{value
function}, where, beginning from an initial point $x_{0}$ at time
$0$
\begin{equation}
V\left(x_{0}\right)=\min_{u_{[0,T]}}\mathbb{E}\left[J\left(x_{0}\right)\right]\label{eq:value-def}
\end{equation}
where we use the shortand $u_{[t,T]}$ to denote the trajectory of $u_t$ over the time interval $t\in[t,T]$, and the expectation is taken over a realizations of the trajectory $x_{\left[0,T\right]}$ with initial condition $x_0$.

The associated Hamilton-Jacobi-Bellman equation, arising from Dynamic Programming arguments \cite{Fleming:2006tl}, is
\begin{equation}
-\partial_{t}V=\min_{u}\left(r+\left(\nabla_{x}V\right)^{T}f+\frac{1}{2}Tr\left(\left(\nabla_{xx}V\right)G\Sigma_{\epsilon}G^{T}\right)\right). \label{eq:HJB-raw}
\end{equation}
As the control effort
enters quadratically into the cost function it is a simple matter
to solve for it analytically by substituting (\ref{eq:cost}) into
(\ref{eq:HJB-raw}) and finding the minimum, yielding:
\begin{equation}
u^{*}=-R^{-1}G^{T}\left(\nabla_{x}V\right).
\label{eq:optimal-u}
\end{equation}

The minimal control, $u^{*}$, may then be substituted into (\ref{eq:HJB-raw})
to yield the following nonlinear, second order partial differential
equation
\begin{equation}
\begin{split}
-\partial_{t}V &=   q  +\left(\nabla_{x}V\right)^{T}f-\frac{1}{2}\left(\nabla_{x}V\right)^{T}GR^{-1}G^{T}\left(\nabla_{x}V\right) \\
 &\phantom{=} +\frac{1}{2} Tr\left(\left(\nabla_{xx}V\right)B\Sigma_{\epsilon}B^{T}\right).
 \end{split}
 \label{eq:HJB-nonlinear}
\end{equation}

The difficulty of solving this PDE is what usually prevents the value function from being directly solved for.
However, it has recently been found \cite{Kappen:2005bn,Todorov:2009je}
that with the assumption that there exists a $\lambda\in\mathbb{R}$
and a control penalty cost $R\in\mathbb{R}^{n\times n}$ satisfying
this equation
\begin{equation}
\lambda G(x)R^{-1}G(x)^{T}=B(x)\Sigma_{\epsilon}B(x)^{T}\triangleq\Sigma_{t}\label{eq:noise-assumption}
\end{equation}
and using the logarithmic transformation, with $\lambda > 0$,
\begin{equation}
V=-\lambda\log\Psi\label{eq:log-transform}
\end{equation}
it is possible, after substitution and simplification, to obtain the
following linear PDE from Equation \eqref{eq:HJB-nonlinear}
\begin{equation}
-\partial_{t}\Psi=-\frac{1}{\lambda}q\Psi+f^{T}\left(\nabla_{x}\Psi\right)+\frac{1}{2}Tr\left(\left(\nabla_{xx}\Psi\right)\Sigma_{t}\right).
\label{eq:HJB-linear}
\end{equation}
This transformation of the value function, which we call here the
\emph{desirability} \cite{Todorov:2009je}, provides an additional, computationally
appealing, method by which to calculate the value function.

\begin{rem}
The condition \eqref{eq:noise-assumption} can roughly be interpreted
as a controllability-type condition: the system controls must span (or counterbalance) the effects of input noise on the system dynamics. A degree of designer input is also given up, as the constraint restricts the design of the control penalty R, requiring
that control effort be highly penalized in subspaces with little noise,
and lightly penalized in those with high noise. Additional discussion
may be found in \cite{Todorov:2009je}.
\end{rem}

The boundary conditions of \eqref{eq:HJB-linear} correspond to the
exit conditions of the optimal control problem. This may correspond
to colliding with an obstacle or goal region, and in the finite horizon
problem there is the added boundary condition of the terminal cost
at $t=T$. These final costs must then be transformed according to
\eqref{eq:log-transform}, producing added boundary conditions to
\eqref{eq:HJB-linear}.

Linearly solvable optimal control is not limited to the finite horizon
setting. Similar analysis can be performed to obtain linear HJB PDEs
for infinite horizon average cost, and first-exit settings, with the
corresponding cost functionals and PDEs shown in Table \ref{tab:control_pdes}. 

\begin{table}

\caption{Linear Desirability PDE for Various Stochastic
Optimal Control Settings. $L(\Psi):=f^{T}\left(\nabla_{x}\Psi\right)+\frac{1}{2}Tr\left(\left(\nabla_{xx}\Psi\right)\Sigma_{t}\right)$}
\label{tab:control_pdes}
\begin{centering}
\begin{tabular}{|c|c|c|}
\hline 
 & Cost Functional & Desirability PDE \\
\hline 
Finite & $\phi_{T}(x_{T})+\int_{0}^{T}r(x_{t},u_{t})dt$ & $\frac{1}{\lambda}q\Psi-\frac{\partial\Psi}{\partial t}=L(\Psi)$\\
\hline 
First-Exit & $\phi_{T_{*}}(x_{T_{*}})+\int_{0}^{T}r(x_{t},u_{t})dt$ & $\frac{1}{\lambda}q\Psi=L(\Psi)$\\
\hline 
Average & $\lim_{T\to\infty}\frac{1}{T}\mathbb{E}\left[\int_{0}^{T}r(x_{t},u_{t})dt\right]$ & $\frac{1}{\lambda}q\Psi-c\Psi=L(\Psi)$\\
\hline
\end{tabular}
\end{centering}
\end{table}

\section{Separated Representations of Tensors}
\label{sec:SR_Tensors}
Traditional numerical techniques to solve PDEs rely on the discretization
of the domain. However, in these schemes the degrees of freedom in the problem grows exponentially with the number of dimensions. While tractable when the number of dimensions is small, in higher dimensions these problems become computationally prohibitive. In \cite{Beylkin:2005kh}, Beylkin and Mohlenkamp proposed to model the solutions to such problems via so-called separated representations, which may be viewed as an extension of the separation of variables technique. Problem data, and the solution, is modeled as a sum of terms, each of which is dependent on individual dimensional variables. Specifically, a function is modeled as
\begin{equation}
f\left(x_{1},\ldots,x_{d}\right)\approx \sum_{l=1}^{r}s_{l}\phi_{1}^{l}(x_{1})\cdots\phi_{d}^{l}(x_{d}).
\label{eq:SR-def}
\end{equation}

The key is that such a representation separates the dependence
of the solution into each component dimension. By then framing operations
to act on single dimensions, it is possible to create algorithms
that need only operate along each dimension independently and thus
scale linearly in $d$. 

However, the complexity of the problem now
grows with $r,$ denoted the \emph{separation rank}. Thus, maintaining a low separation rank becomes paramount for any practical algorithm. Unfortunately, many operations inherently increase the separation rank, including vector addition and matrix-vector multiplication. This unbounded growth is mitigated by reducing the separation rank at each step of an algorithm in an attempt to continually maintain low rank approximations. Unfortunately, there are often no guarantees that a given function, or solution to a PDE, will have low separation rank and situations may arise where it is impossible to lower the rank while maintaining the desired accuracy.

Here we simply provide an introduction to the concept of a separated representation and discuss them in a manner tailored to our use. As such, we direct the reader to \cite{Beylkin:2005kh} for a complete treatment. 

A vector $F$ in dimension $d$ is a discrete representation of a
function $f$ on a rectangular domain, $\boldsymbol{F}=F(j_{1},\ldots,j_{d})$
where $j_{i}=1,\ldots,M_{i}$ are the indices along each dimension. A linear operator $\mathcal{A}$ in dimension $d$ is a linear map $\mathcal{A}:S\to S$ where $S$ is the
space of functions in dimension $d$. A matrix $\mathbb{A}$ in dimension
$d$ is a discrete representation of a linear operator in dimension
$d.$
\begin{defn}
For a given $\epsilon,$ we represent a vector $\boldsymbol{F}=F(j_{1},j_{2},\ldots,j_{d})$
in dimension $d$ as 
\[
\boldsymbol{F} \approx \sum_{l=1}^{r_{\textbf{F}}}s_{l}\bigotimes_{i=1}^d\boldsymbol{F}_{i}^{l}
\]
where $\bigotimes$ denotes the tensor product and $\boldsymbol{F}_{i}^{l}$
are traditional vectors in $\mathbb{R}^{M_i}$ with entries $F_{i}^{l}(j_{i})$ and
unit norm. For this to be an $\epsilon$ accurate representation we require that
\[
\left\Vert \boldsymbol{F}-\sum_{l=1}^{r_{\textbf{F}}}s_{l}\bigotimes_{i=1}^d\boldsymbol{F}_{i}^{l}\right\Vert \le\epsilon.
\]
The integer $r$ is known as the \textbf{separation rank}.
\end{defn}
The matrix definition is analogous, with the matrices ${\mathbb{A}_{i}^{l} \in \mathbb{R}^{M_i \times M_i}}$
in lieu of $\boldsymbol{F}_{i}^{l}=F_{i}^{l}(j_{i})$. Matrix multiplication is then performed as 
\begin{equation}
\mathbb{A}F=\sum_{m=1}^{r_{\mathbb{A}}}\sum_{l=1}^{r_{F}}s_{m}^{\mathbb{A}}s_{l}^{F}\left(\mathbb{A}_{1}^{m}F_{1}^{l}\right)\otimes\cdots\otimes\mathbb{A}_{d}^{m}F_{d}^{l}.
\end{equation}
Since matrix operations in this formulation reduce to individual operations
along each dimension, as the dimensionality of the problem increases the complexity of these operations scales linearly, \emph{e.g.}, if we let $M_i$ = $M$ for all $i$ a matrix vector multiplication costs $O(r_A r_F d M^2)$. 

\subsection{Alternating Least Squares}
Any scheme that uses these separated representations will become computationally prohibitive if the separation ranks are allowed to grow unchecked. For example, in the matrix vector multiplication the separation rank of the output grows by a factor of $r_A r_F,$ so even performing the most basic of operations may have a large impact on the separation rank. It is therefore necessary to periodically find lower rank approximations to various operators and vectors. If the assumption is that the discrete versions of the functions being represented have low separation rank, then any increase in the separation rank is only an artifact of the way operations are performed in these tensor representations, and we expect is possible to produce an accurate representation of the resultant tensor that has a reduced separation rank.

We provide a high level overview of the alternating least squares (ALS) algorithm, with the reader directed to \cite{Beylkin:2005kh} for details. This algorithm allows us to look for separated representations with lower separation rank and solve systems. A recently proposed variant relying on a randomized interpolative decomposition is presented in \cite{Biagioni:2013wh} and may be used as a precursor to ALS. 

Given a vector with separated representation $\textbf{G}$, and an operator $\mathbb{A}$, ALS tries to either minimize $\left\Vert \textbf{F} - \textbf{G}\right\Vert$ or, with slight modification, $\left\Vert \mathbb{A} \textbf{F} - \textbf{G}\right\Vert$ subject to $\textbf{F}$ having fixed rank $r_{F}.$ The algorithm sweeps through the coordinate directions, effectively
performing block-gradient descent. For a fixed separation rank of $\textbf{F}$ this process may be repeated until the algorithm has either achieved the desired accuracy, or has stagnated. If the algorithm has stagnated, and the representation error is not small enough, \emph{e.g.}, $\left\Vert \mathbb{A} \textbf{F}-\textbf{G}\right\Vert \ge\epsilon,$ a random rank-one tensor is added to $\textbf{F}$, and the ALS routine continues.

At an individual step in this iterative algorithm, all dimensions of the tensor $\textbf{F}$ are held constant save one dimension $k,$ in which case the least-squares problem becomes linear in $\textbf{F}_{k}^{l}$ for $l=1,\ldots,r_G$. The resulting least squares problem is solved via derivation of the normal equations, yielding a linear set of equations. These were originally given in \cite{Beylkin:2005kh}, but have been rederived in \cite{Sun2014} in a more compact format as

\begin{equation}
\left(\begin{array}{ccc}
M_{1,1} & \cdots & M_{1,r_{G}}\\
\vdots & \ddots & \vdots\\
M_{r_{G},1} & \cdots & M_{r_{G},r_{G}}
\end{array}\right)\left(\begin{array}{c}
F_{k}^{1}\\
\vdots\\
F_{k}^{r_{G}}
\end{array}\right)=\left(\begin{array}{c}
N_{1}\\
\vdots\\
N_{r_{G}}
\end{array}\right)
\end{equation}
where the components of the normal equations are given by
\begin{equation}
M_{i,j}=\sum_{i_{A}=1}^{r_{\mathbb{A}}}\sum_{j_{A}=1}^{r_{\mathbb{A}}}\left(A_{k}^{j_{A}}\right)^{T}A_{k}^{i_{A}}\prod_{d\neq k}\left\langle A_{d}^{i_{A}}F_{d}^{j},A_{d}^{j_{A}}F_{d}^{i}\right\rangle 
\end{equation}
\begin{equation}
N_{i}=\sum_{i_{A}=1}^{r_{\mathbb{A}}}\sum_{i_{G}=1}^{r_{G}}\left(A_{k}^{i_{A}}\right)^{T}G_{k}^{i_{G}}\prod_{d\neq k}\left\langle A_{d}^{i_{A}}F_{d}^{i},G_{d}^{i_{G}}\right\rangle 
\end{equation}
with the solution vector $
\textbf{F}=\sum_{i=1}^{r_{\textbf{F}}}\otimes_{d=1}^{n}F_{d}^{i}
$, the vector $
\textbf{G}=\sum_{i=1}^{r_{\textbf{G}}}\otimes_{d=1}^{n}G_{d}^{i}
$, and the operator $
\mathbb{A}=\sum_{i=1}^{r_{\mathbb{A}}}\otimes_{d=1}^{n}A_{d}^{i}
$, for $n$ the dimension of the system. When the operator $\mathbb{A}$ is the identity, the problem $\left\Vert \textbf{F}-\textbf{G}\right\Vert $ has additional structure that may be leveraged \cite{Beylkin:2005kh}. The solution to these linear equations provides a new $\textbf{F}$, with different components in the $k$ dimension, that has a smaller residual error. The algorithm continues by optimizing the components in the next, $k+1,$ dimension.

While the core ALS algorithm, with the identity operator, costs $\mathcal{O}\left(dM+dr_{\mathbf{F}}^{3}\right)$ per iteration, its use to solve a linear system costs $\mathcal{O}\left(dM^{3}+r_{\mathbb{A}}^{3}M^{3}\right)$ per iteration, where $d$ is the underlying dimensionality of the system, $M$ is the maximal number of mesh nodes along each dimension, and $r_\mathbb{A}, (r_\mathbf{F})$  is the rank of the operator $\mathbb{A}$ (vector $\mathbf{F}$). See \cite{Beylkin:2005kh} for a more comprehensive list of algorithms that may be used with operators and vectors in separated representations. 

\section{Separated Solution to the HJB}
\label{sec:SR_HJB}
We make the modeling assumption that the problem data of (\ref{eq:dynamics}) can be accurately represented, or approximated, with a low rank separated  representation.
\begin{equation}
f_{i}(t,x)=\sum_{l=1}^{r_{f_{i}}}\bigotimes_{k=1}^{d}\left(f_{i}\right)_{d}^{l}
\end{equation}
where $r_{f_i}$ is assumed to be small.

There is then the need to approximate the relevant operators present in (\ref{eq:dynamics}), specifically the gradient and Hessian, in a low rank representation. A number of options exist, with varying levels of complexity in the analysis and accuracy, ranging from simple finite difference schemes to spectral
differentiation techniques \cite{Trefethen:2000wd}. Specifically, the gradient along dimension $k$ is simply
\[
\nabla_{k}=I_{1}\otimes\cdots\otimes I_{k-1} \otimes \nabla \otimes I_{k+1} \otimes\cdots\otimes I_{d}
\]
while the Hessian has entries $\nabla_{k,j}=\nabla_{k}\cdot\nabla_{j}$,
and the estimates of the derivative along an individual coordinate
are simply a suitably high order finite difference scheme in one dimension.
Thus, the directional gradient and second order terms may simply be constructed out of rank one representations. For example, using of sums of these rank one terms yields a representation for the Laplacian that has separation rank $d.$ However, such a representation may be not have minimal separation rank for a given accuracy. Other constructions specifically targeting the separated representation exist \cite{Beylkin:2005kh}, for example a Laplacian approximation may be made with separation rank two, rather than requiring a full rank-$d$ sum of second order terms. Choices about how to approximate the operators may lead to some variation in the separation rank of the solutions since different strategies introduce different types of discretization error. 

\subsection{Separation Rank of the HJB}
Determining the separation rank of the HJB operator is straightforward. Denote the rank of a vector or operator $X$ as $r_X$. Recalling (\ref{eq:HJB-linear}) and neglecting the time dependent component, the operator consists of three additive terms. The state-cost term $q \Psi$ is a diagonal operator along each dimension, and thus contributes $r_{q}$. The second, advection term is an inner product between the dynamics $f$ and the gradient, resulting in the multiplication of each element
$f_{i}$ by a rank one operator, and then their summation. The contribution
from this component results in separation rank $\sum_{k=1}^{d}r_{f_{i}}$. 
Finally, the second-order term requires the construction of $\Sigma_{t}$
in (\ref{eq:noise-assumption}). Here the growth in the separation rank may be significant, due to the multiplicative contribution of $G$. However, given diagonal cost matrix $R$ or noise covariance $\Sigma_\epsilon$ the number of terms may collapse significantly. The separation rank of the HJB operator is simply the sum of these three terms' rank.

The result is that the separation rank for individual problems may
vary over a wide range, depending on the problem data. However, in
many problems of interest it remains low. For even apparently complex systems,
complexity typically manifests as nonlinear multiplicative terms in
the dynamics. This form of complexity effectively adds no cost in terms of separation rank, and it is instead the number of additive terms that are of concern, which
is typically small. Furthermore, in many applications the control and
noise matrices typically contain constant terms, corresponding to tensors
of separation rank one. Finally, for systems where a high separation
rank accumulates, it remains possible to search for low rank structure
by performing ALS on the operator before attempting to solve the linear
system.

\subsection{Representation of Interior Boundary Conditions}
\label{sec:Interior_BCs}




Optimal control applications impose irregular boundary conditions on many problems of interest. For example, stabilization to the origin corresponds to a zero-cost point-boundary at the origin. Obstacles, or sub-task goals in temporal problems (\cite{Esfahani:2012wf}, \cite{Horowitz:2014we}), become boundary conditions as well, i.e. a set of Dirichlet boundary conditions within the domain.

We impose \emph{essential} boundary conditions by setting the value of nodes to some desired value via linear equalities within the domain. Although in other settings it is desirable to remove the degrees of freedom from within the boundaries to save computational effort, in our context maintaining the grid form is a far greater concern. Specifically, we impose Dirichlet boundary conditions only on regions composed of hyper-cubes in the domain, allowing us to modify the domain with only a modest increase in the separation rank of the operator. We first eliminate the operators effect on this hypercube of the domain via a single subtraction, and then replace it with the identity operator via a single addition. The resulting operator $\bar{\mathbb{A}}$ with the desired boundary conditions has rank $r_{\bar{\mathbb{A}}} = d+2r_{\mathbb{A}}$. 

\section{Implementation Details and Examples}
\label{sec:Examples}

In the following examples, we approximate first and second order derivatives using eighth order finite differences, with the number of mesh points along each dimension varying between $n_g=100$ and $n_g=201$. The result are tensors that would typically not fit in the memory of even the largest modern super computers if expressed naively without the use of the separated representation. In each case we modeled the problem as first-exit (see Table \ref{tab:control_pdes}). In all cases the noise was assumed to enter the dynamics in the same manner as the control, with $G(x)=B(x)$ in (\ref{eq:dynamics}).

The operator is constructed as described in Section \ref{sec:SR_HJB}. The operator and boundary conditions are compressed independently using Alternating Least Squares with the linear system set to identity. With this low-rank representation, the problem is then solved using Alternating Least Squares for the HJB system. We employed the Matlab Tensor Toolbox \cite{TTB_Software, TTB_CPOPT}, for storage and manipulation of tensor objects.

The problems were solved on a quad-core 2.3Ghz i7-equipped laptop. We denote $\bar{u}$, $\bar{x}$ as the vector of system control inputs and states for each example.

\subsection{Inverted Pendulum}

In \cite{Osinga:2006bk} the geometry of optimal control for the inverted pendulum on a cart was investigated in detail. In particular, they produce the value function for the inverted pendulum when actuated directly at the base
\begin{eqnarray}
\dot{x}_{1} & = & x_{2} \\
\dot{x}_{2} & = & \frac{\frac{g}{l}\sin(x_{1})-\frac{1}{2}m_{r}x_{2}^{2}\sin(2x_{1})-\frac{m_{r}}{ml}\cos(x_{1})u}{\frac{4}{3}-m_{r}\cos^{2}(x_{1})} \nonumber
\end{eqnarray}
and the cost function is $q(x) = 0.1 x_1^2 + 0.05 x_2^2 + 0.01 u^2$. This problem has periodic boundary conditions along the $x_1$ dimension, and we placed a Dirichlet boundary condition of $\phi(x_1,\pm 11)=10$, i.e. a high penalty for exceeding the maximal angular velocity of $\dot{\theta}>11$ rad/s. An exit interior boundary was placed at the origin, with Dirichlet boundary conditions corresponding to unity desirability. We chose $n_g=201$ discretization points in each dimension.

\begin{figure}
\begin{center}
\includegraphics[scale=0.45, trim = 0mm 10mm 0mm 0mm]{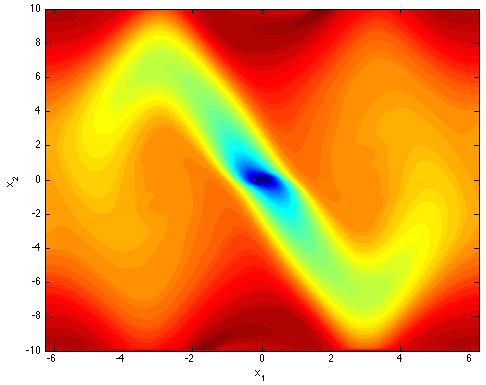}
\end{center}
\caption{Cost to go for the inverted pendulum. The effects of the noise may be seen in the smoothing of the value function, in comparison to the deterministic case seen in \cite{Osinga:2006bk}.
\label{fig:inv_pend_c2go}}
\end{figure}

\begin{figure}
\begin{centering}
\begin{tabular}{c c c c c c c}
$x_1$ &
\includegraphics[scale=0.07, trim = 50mm 30mm 50mm 30mm]{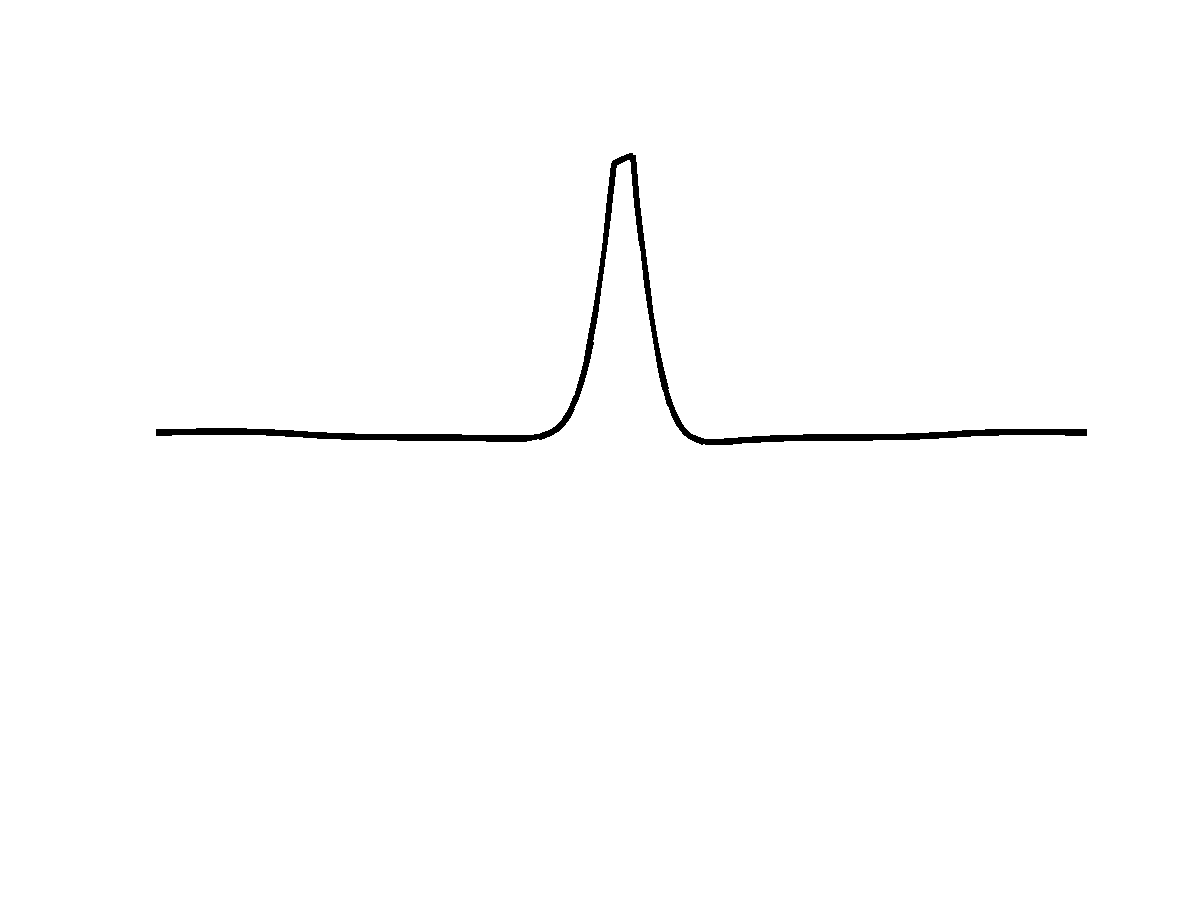} &
\includegraphics[scale=0.07, trim = 50mm 30mm 50mm 30mm]{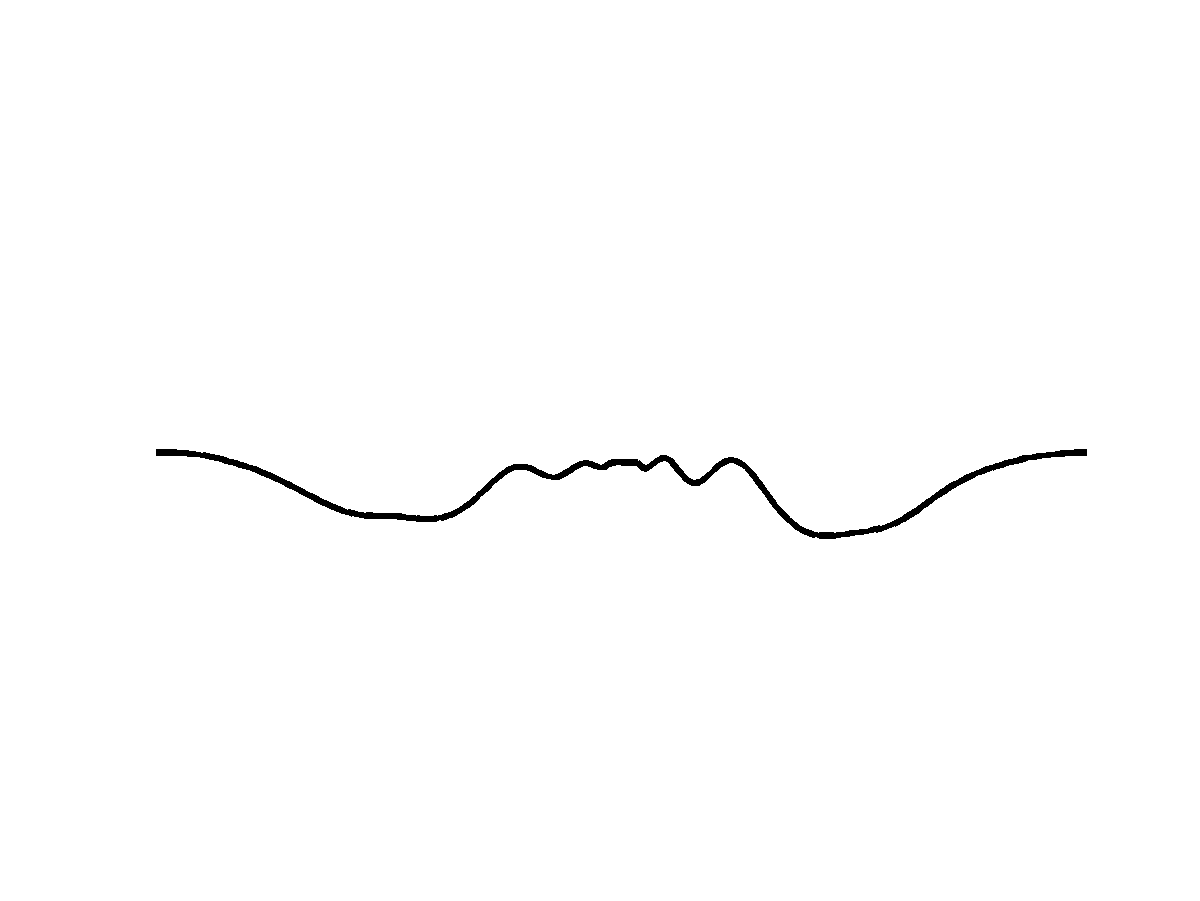} &
\includegraphics[scale=0.07, trim = 50mm 30mm 50mm 30mm]{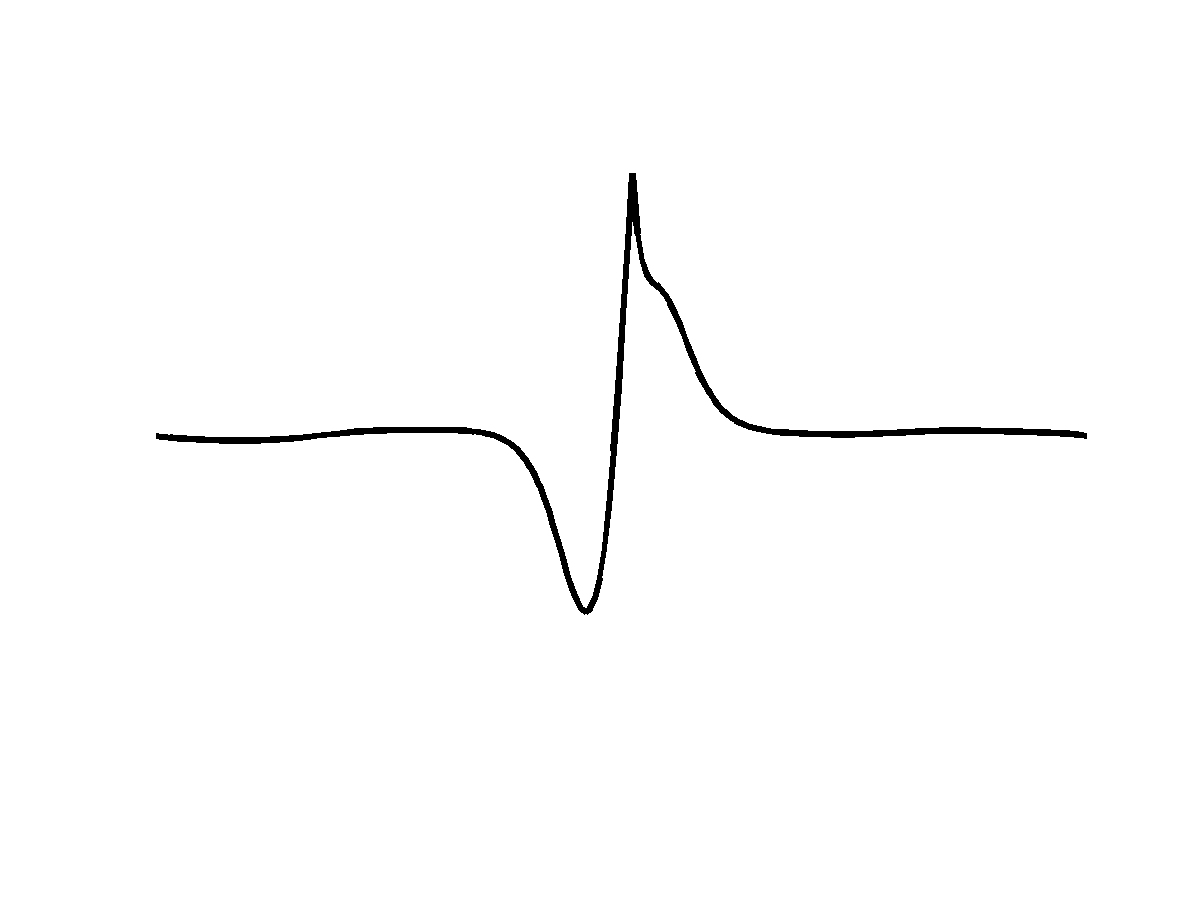} &
\includegraphics[scale=0.07, trim = 50mm 30mm 50mm 30mm]{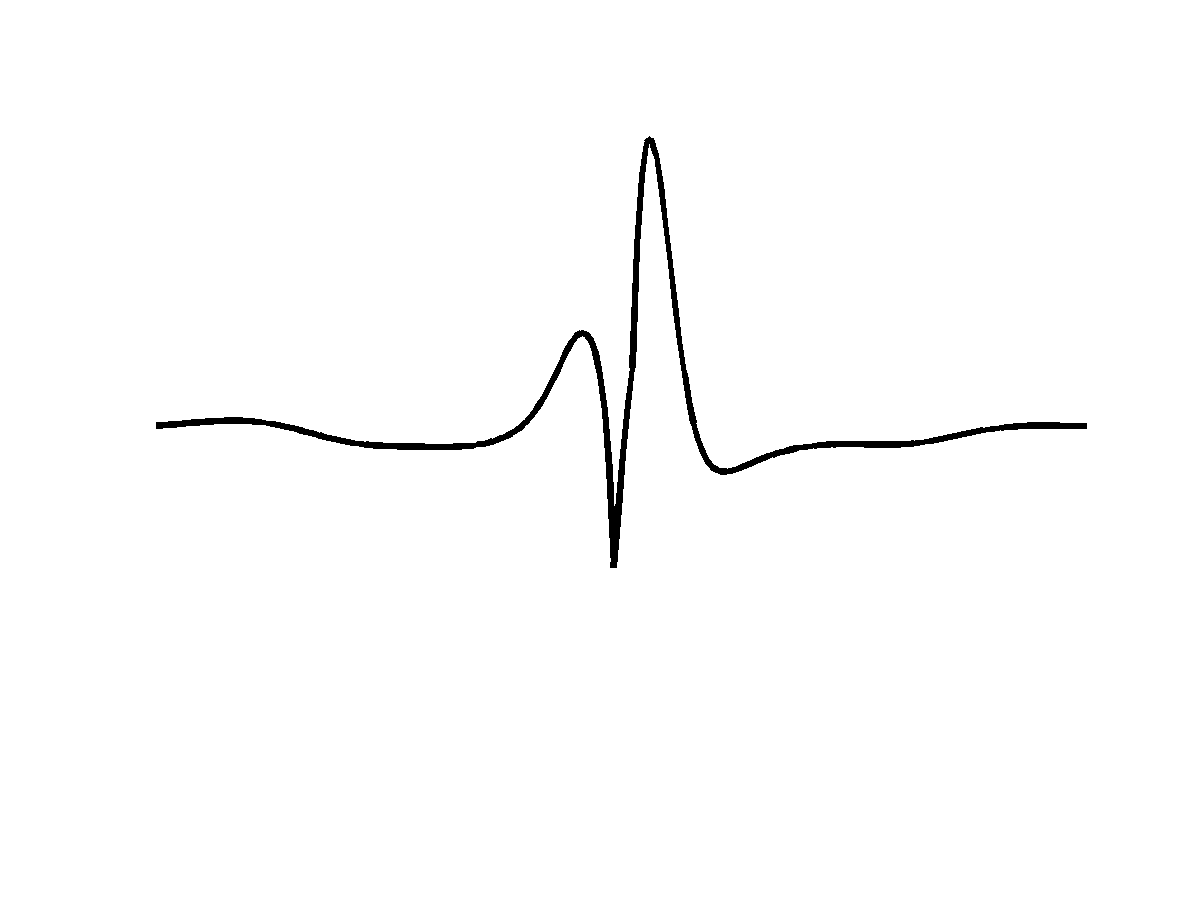} &
\includegraphics[scale=0.07, trim = 50mm 30mm 50mm 30mm]{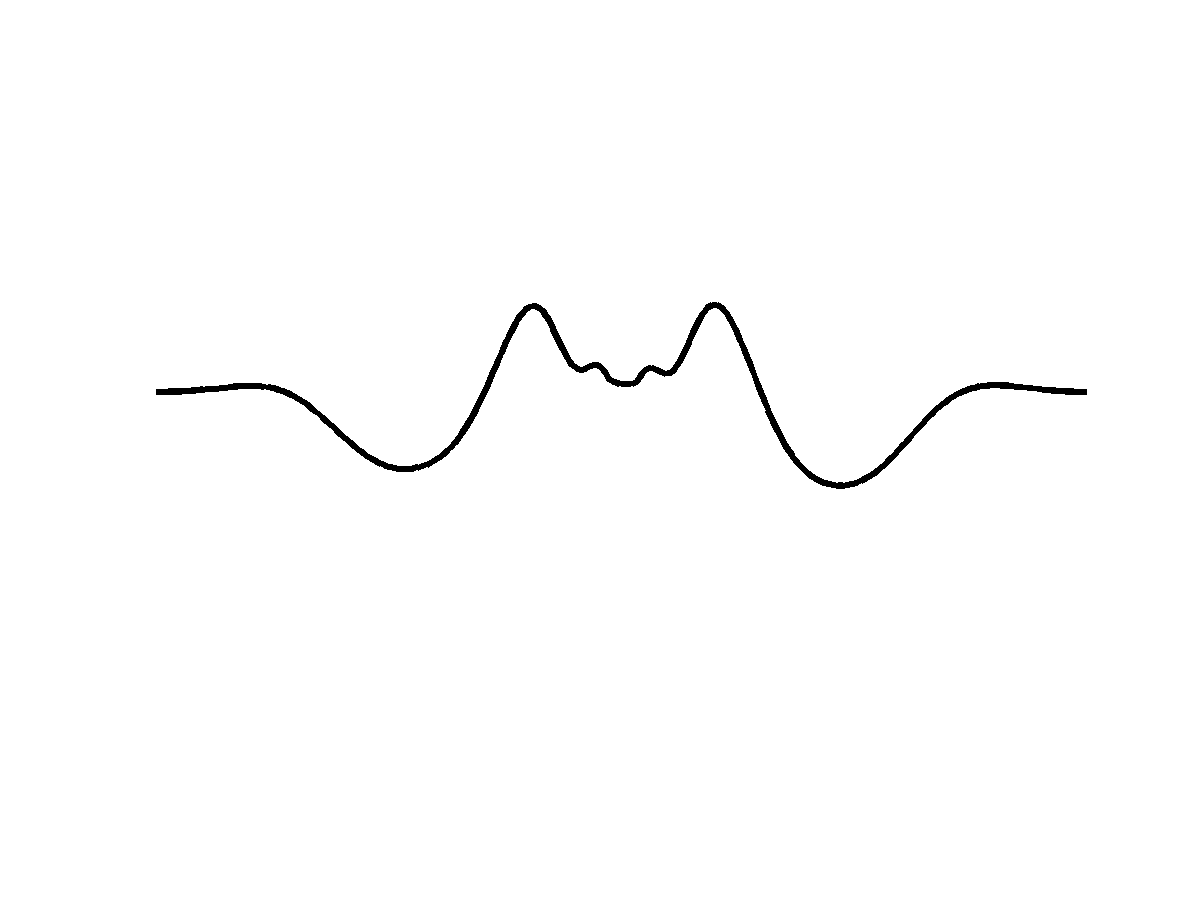} \\

$x_2$ &
\includegraphics[scale=0.07, trim = 30mm 50mm 20mm 30mm]{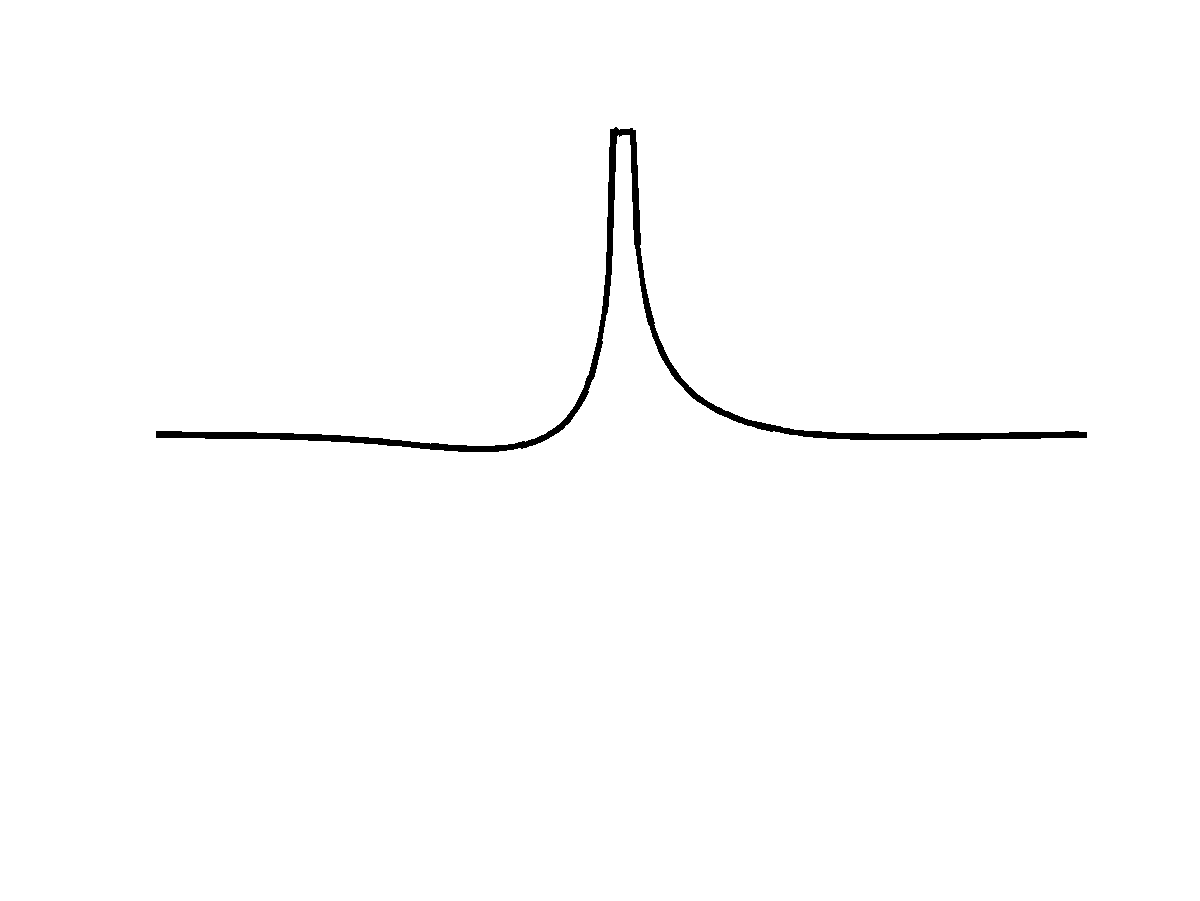} &
\includegraphics[scale=0.07, trim = 30mm 50mm 20mm 30mm]{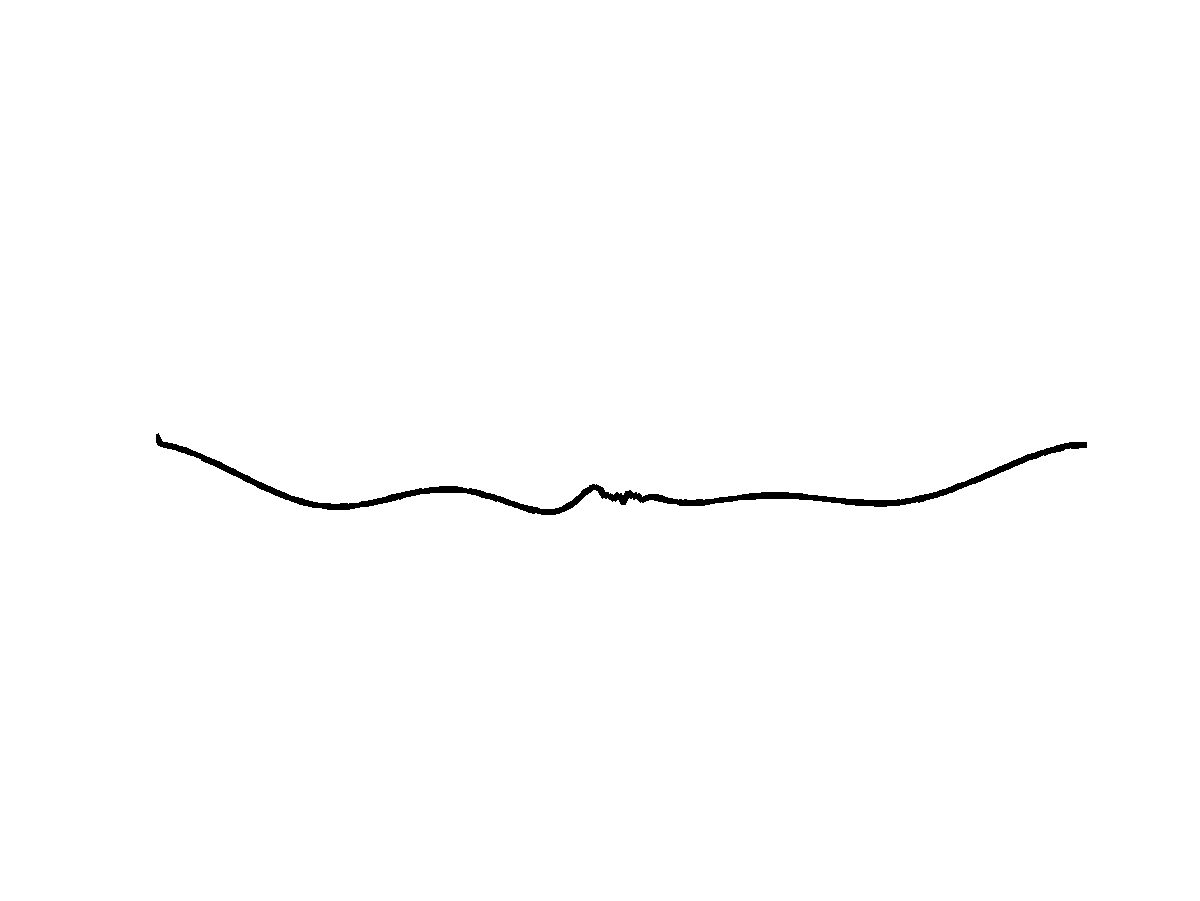} &
\includegraphics[scale=0.07, trim = 30mm 50mm 20mm 30mm]{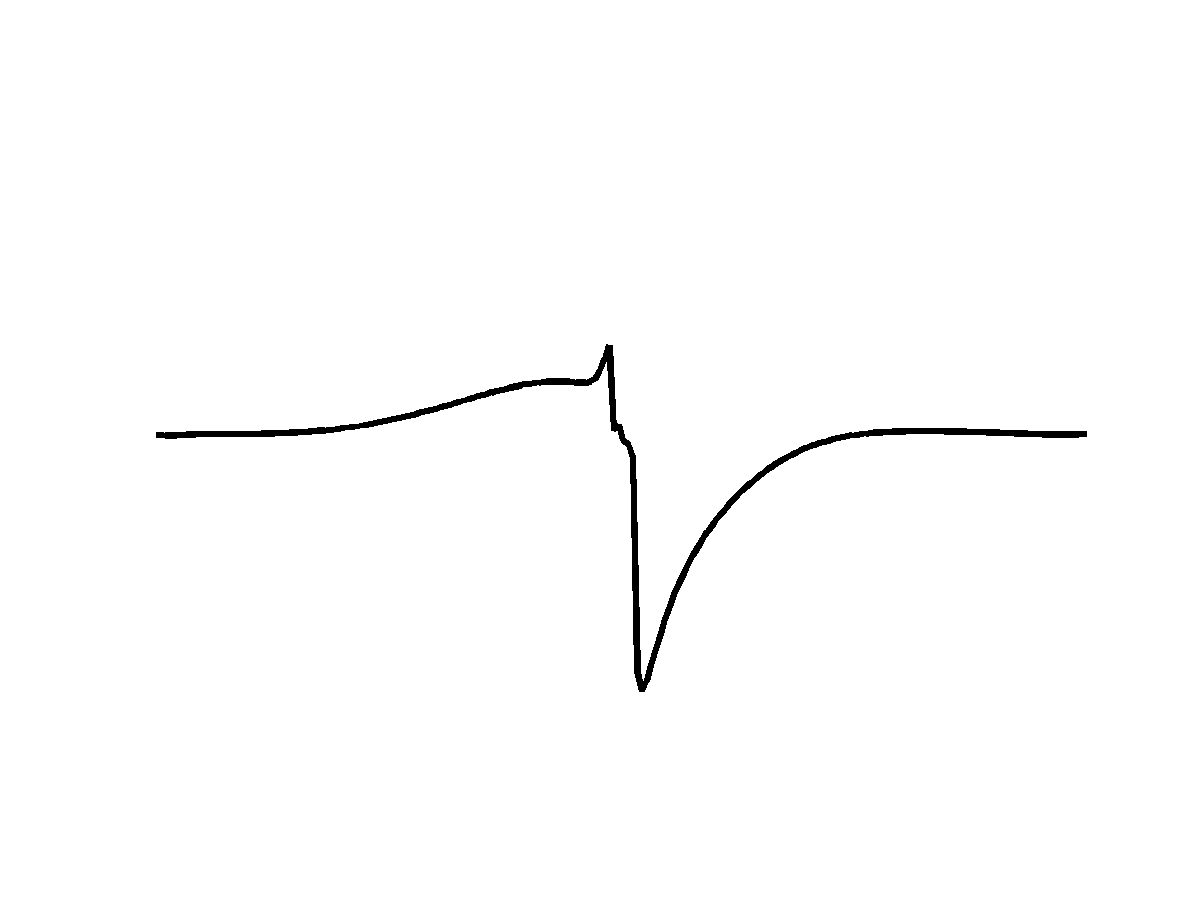} &
\includegraphics[scale=0.07, trim = 30mm 50mm 20mm 30mm]{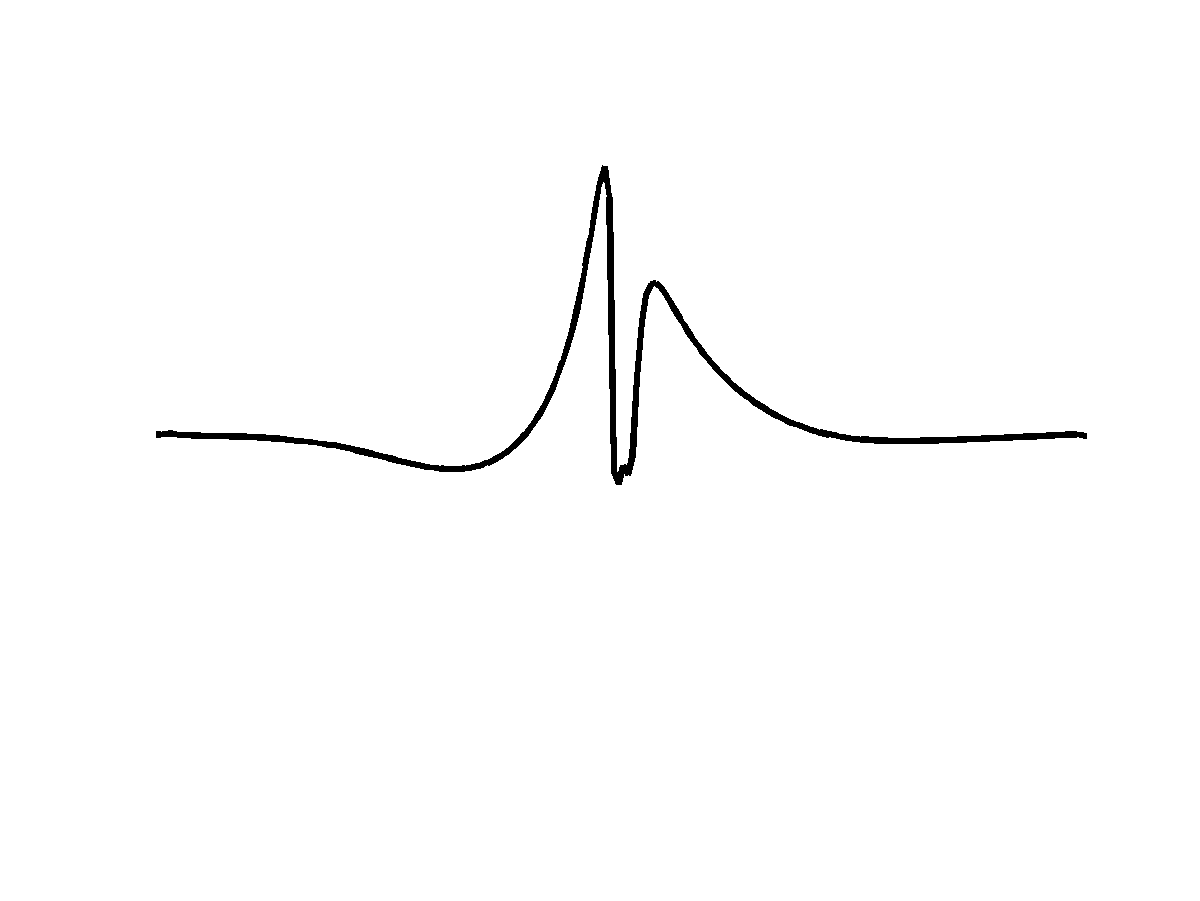} &
\includegraphics[scale=0.07, trim = 30mm 50mm 20mm 30mm]{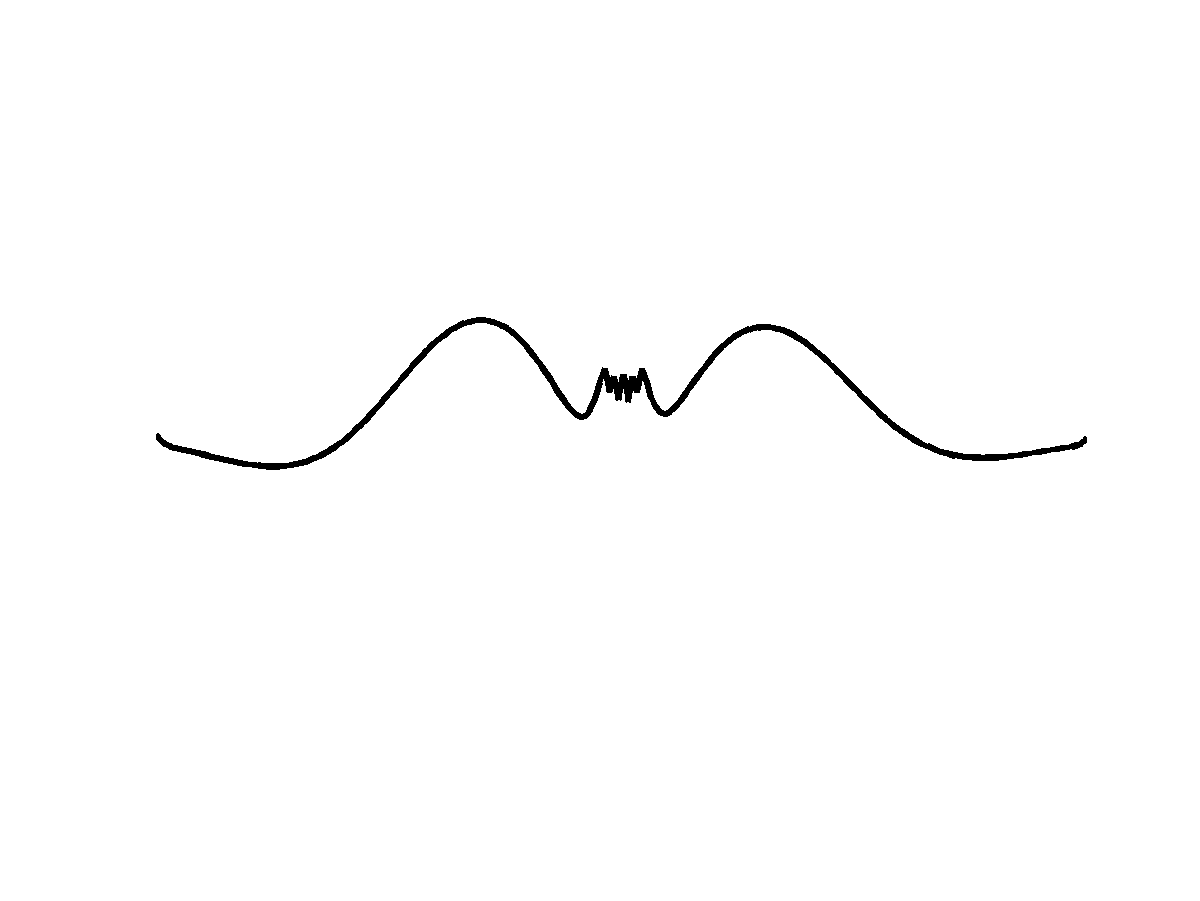}
\end{tabular}
\end{centering}
\caption{Five principal basis functions for the inverted pendulum along the $x_1$, $x_2$ dimensions.
\label{fig:inv_pend_bases}}
\end{figure}

The value function obtained by inverting the transformation (\ref{eq:log-transform}) to the solution is shown in Figure \ref{fig:inv_pend_c2go}. The process took approximately ten minutes, achieving error $e=5.22\cdot 10^{-5}$ with a basis of $r_\Psi=20$ rank one tensors. The five principal basis functions along each dimension are shown in Figure \ref{fig:inv_pend_bases}.

\subsection{VTOL Aircraft}
\begin{figure}
\begin{centering}
\includegraphics[scale=0.48, trim = 10mm 10mm 10mm 0mm]{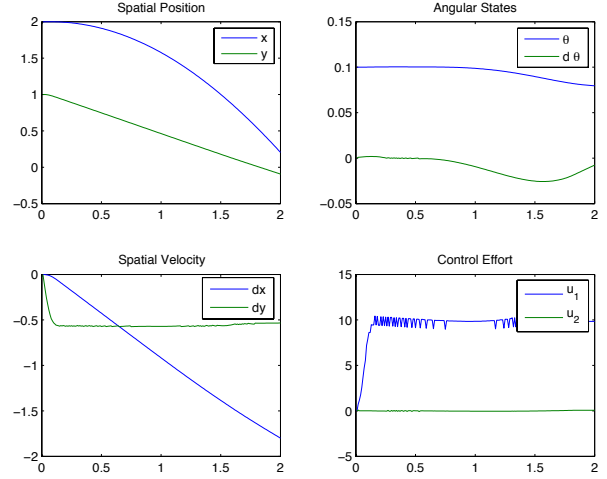}
\end{centering}
\caption{Sample trajectory when executing desirability for the VTOL aircraft.
\label{fig:vtol_simulation}}
\end{figure}

\begin{figure}
\begin{centering}
\includegraphics[scale=0.23, trim = 10mm 10mm 10mm 10mm]{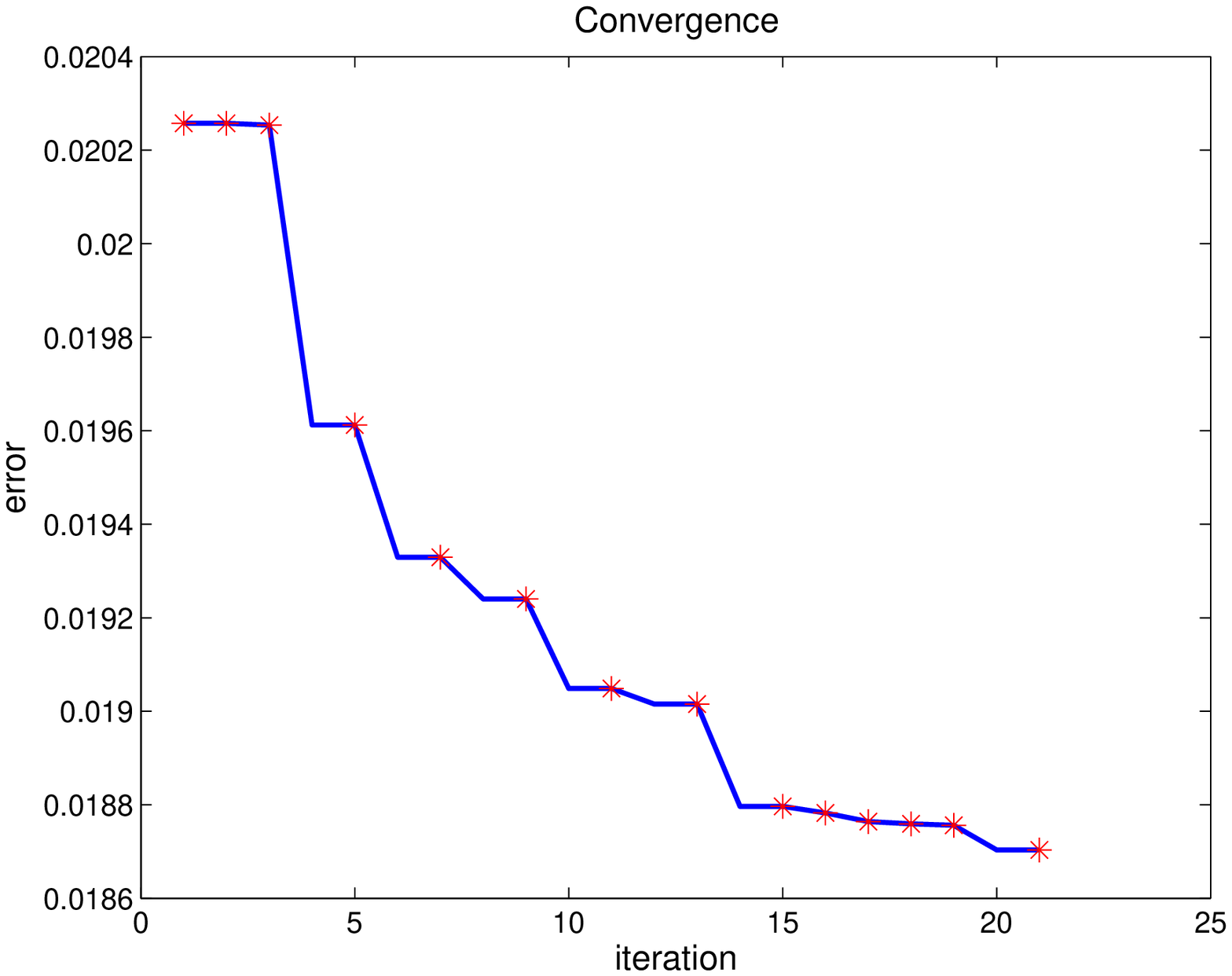}
\includegraphics[scale=0.23, trim = 10mm 10mm 10mm 10mm]{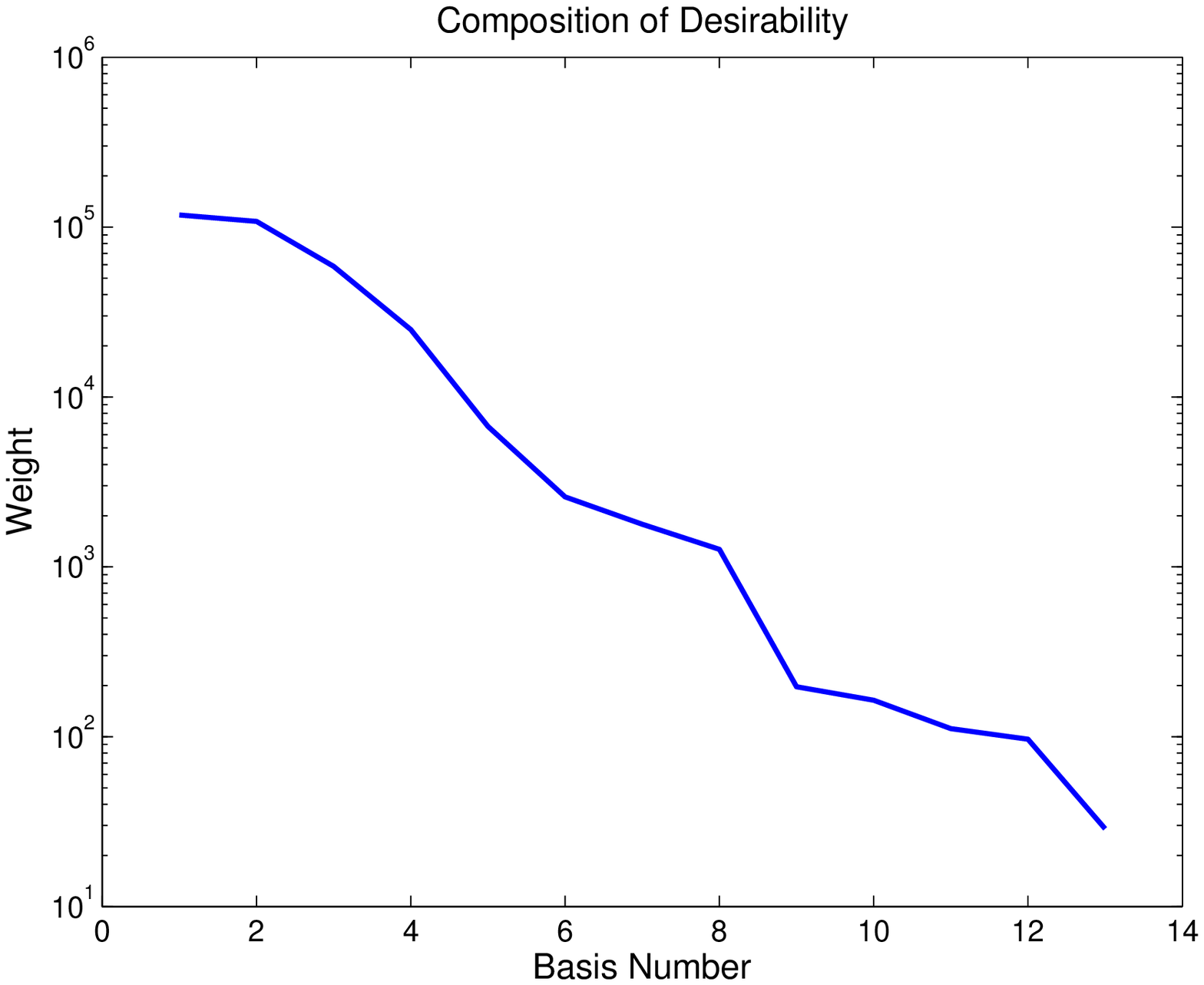}
\end{centering}
\caption{Convergence and weighting for the VTOL solution. The red markers indicate at which iterations the ALS algorithm enriched the solution by adding a basis element.
\label{fig:vtol_convergence}}
\end{figure}

Next, we consider a Vertical Takeoff and Landing aircraft (also known as the Harrier Jet). We examine a planar cross section of the translational state, that is the jet's $(x,y)$ location where $y$ is in the vertical direction. The system is characterized by second order dynamics with gravitational drift and trigonometric inputs, giving rise to a sixth dimensional nonlinear system. Specifically, the equations governing the system are given in \cite{hauser1992nonlinear} as
\begin{eqnarray*}
\ddot{x} & = & -u \sin(\theta) + \epsilon \, \tau \cos(\theta) \\
\ddot{y} & = & u \cos(\theta) + \epsilon \, \tau \sin(\theta) - g \\
\ddot{\theta} & = & \tau,
\end{eqnarray*}
where $\epsilon=0.01$ in our example. The cost function chosen was $r=u^2$, and $q(x,y,\theta,\ldots)=1.0$ on the domain $x\in [-4,4], y\in[0,2], \dot{x}\in[-8,8], \dot{y}\in[-1,1], \dot{\theta}\in[-5,5]$, with $\theta$ periodic on $[-\pi,\pi]$. All boundaries were set to have boundary conditions $\Psi \mid_{\partial \Omega}=0$, save $y=0$, which had condition $\Psi \mid \partial \Omega=1-s^2$ for each coordinate direction $s$, placing a target of landing with zero velocities. Discretization $n_g=100$ were used along each dimension. We limited the solver to twenty iterations, which required approximately five minutes. We omit the resulting basis functions due to space restrictions. We also show the error and basis function weighting in Figure \ref{fig:vtol_convergence}.
 A sample trajectory when executing the policy in closed loop in Figure \ref{fig:vtol_simulation}.

\subsection{Quadcopter}

The next example is in the stabilization of a quadcopter. The derivation of the dynamics may be found in \cite{GarciaCarrillo:2013du}, and results in a system of order twelve with highly nonlinear dynamics.
\begin{eqnarray*}
m\ddot{x} & = & u\left(\sin\phi\sin\psi+\cos\phi\cos\psi\sin\theta\right)\\
m\ddot{y} & = & u\left(\cos\phi\sin\theta\sin\psi-\cos\psi\sin\phi\right)\\
m\ddot{z} & = & u\cos\theta\cos\phi-mg\\
\ddot{\psi} & = & \tilde{\tau}_{\psi}\\
\ddot{\theta} & = & \tilde{\tau}_{\theta}\\
\ddot{\phi} & = & \tilde{\tau}_{\phi}
\end{eqnarray*}
where $\eta=\left(x,y,z\right)$ are in the horizontal and vertical
plane, respectively, while $\tilde{\tau}=\left(\tilde{\tau}_{\psi},\tilde{\tau}_{\theta},\tilde{\tau}_{\phi}\right)$
are the yaw, pitch, and roll moments. For simplicity, we assume we have direct actuation control over $\tilde{\tau}$. We solve the problem with $r=\| \bar{u} \|$ and $q(\bar{x})=2$. Similar to the VTOL example, we penalize all boundaries, save $x=1$, where a quadratic along the boundary in each dimension induces the system to exit with small velocity in all dimensions. Discretization $n_g=100$ was again used along each dimension.

In this instance $f(x)\equiv0$ for all but the $z-$acceleration, which has separation rank one, and $G(x)$ has separation rank two for only the first three coordinate dimensions.
Due to the matching condition (\ref{eq:noise-assumption}), we model
noise as entering the system as entering the same subspace as the
control input, with $B(x)\triangleq G(x)$. The formation of the partial differential operator requires $r_{\mathbb{A}}=56$, but the ALS algorithm is able to compress this to $r_{\tilde{\mathbb{A}}} = 24$ with a relative error of $10^{-4}$ in approximately two minutes, indicating there exist a great deal of underlying structure that the system is able to exploit.

\begin{figure}
\begin{centering}
\begin{tabular}{c c c c}
\includegraphics[scale=0.12, trim = 5mm 5mm 15mm 5mm]{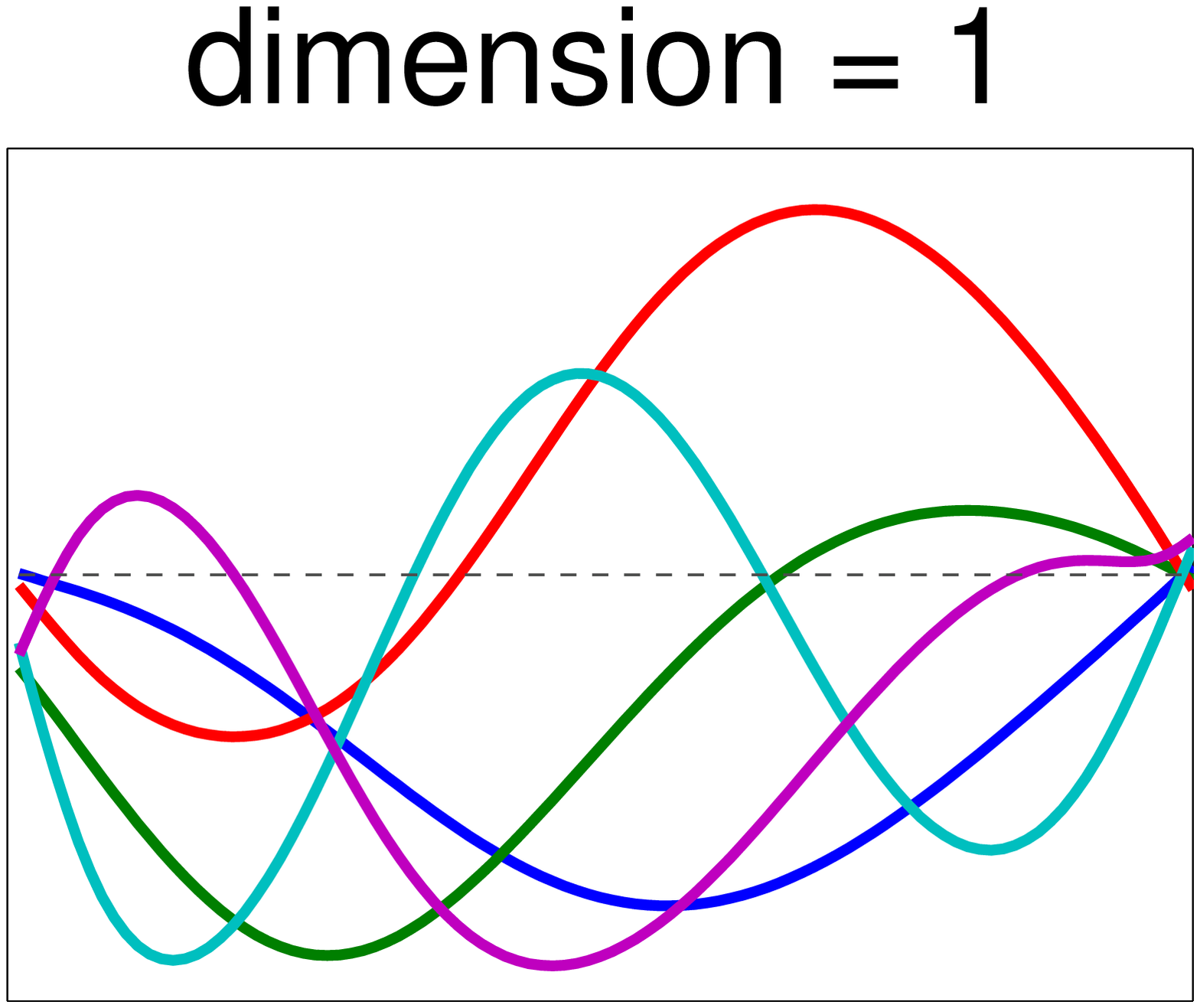} &
\includegraphics[scale=0.12, trim = 5mm 5mm 15mm 5mm]{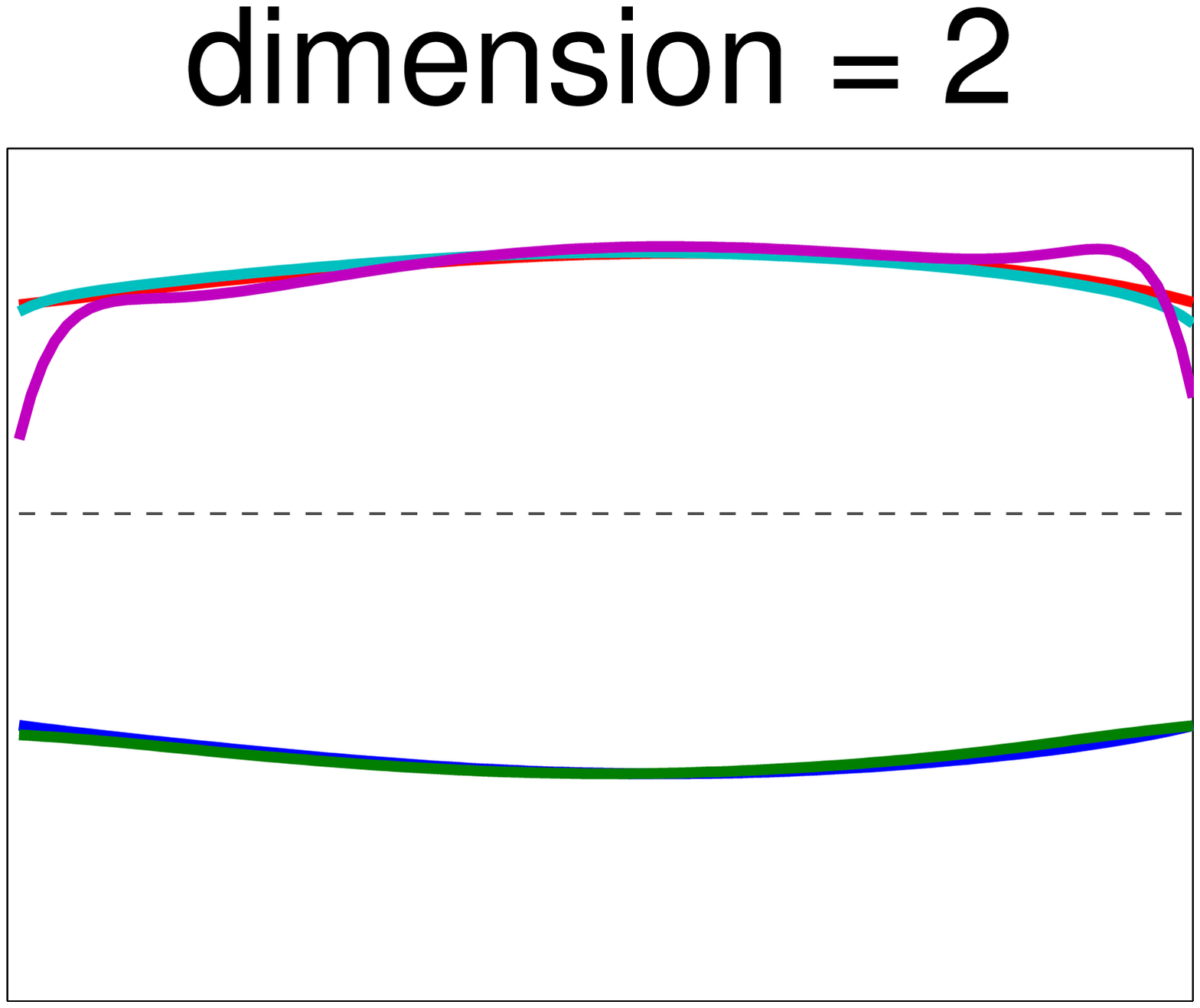} &
\includegraphics[scale=0.12, trim = 5mm 5mm 15mm 5mm]{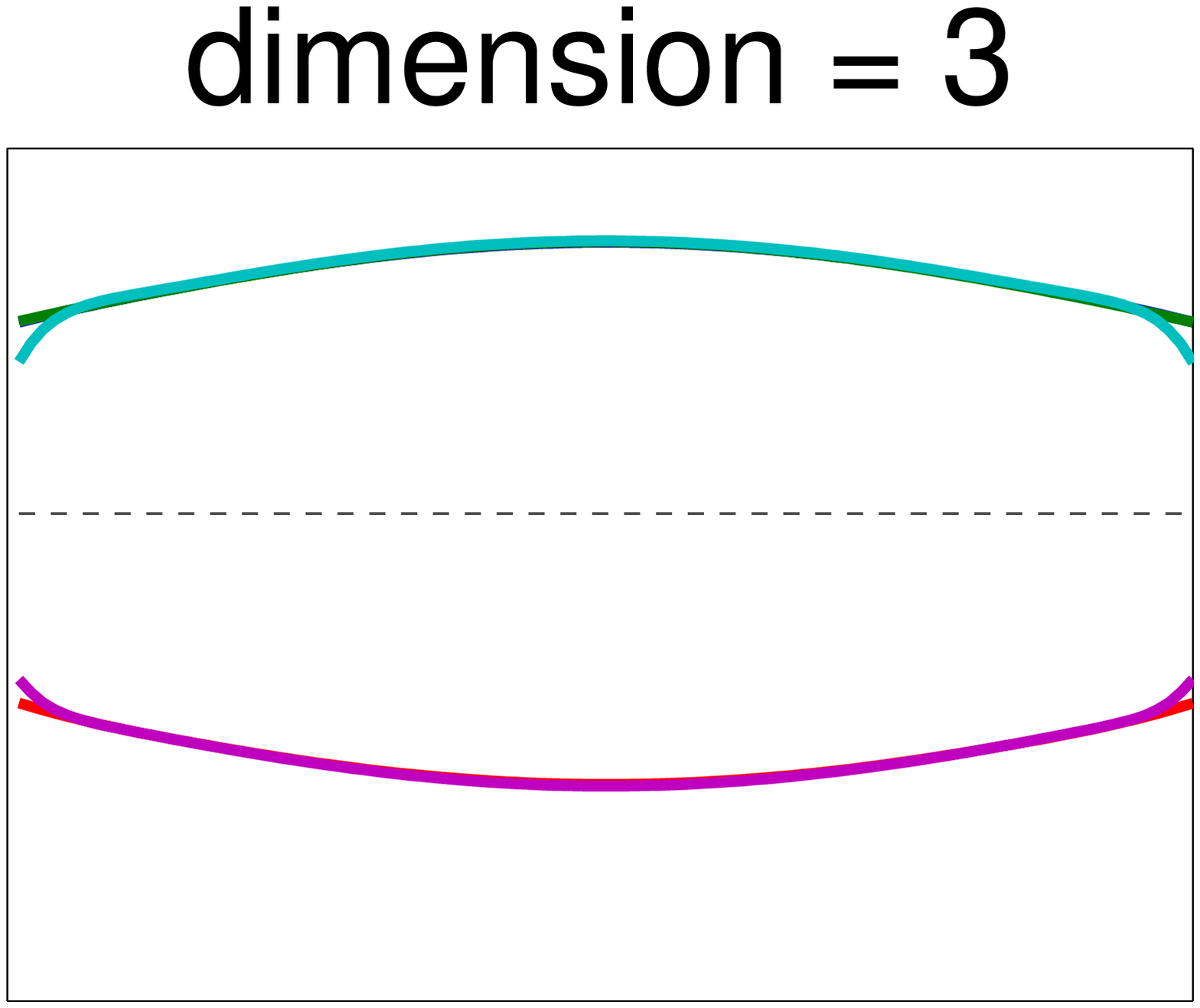} &
\includegraphics[scale=0.12, trim = 5mm 5mm 15mm 5mm]{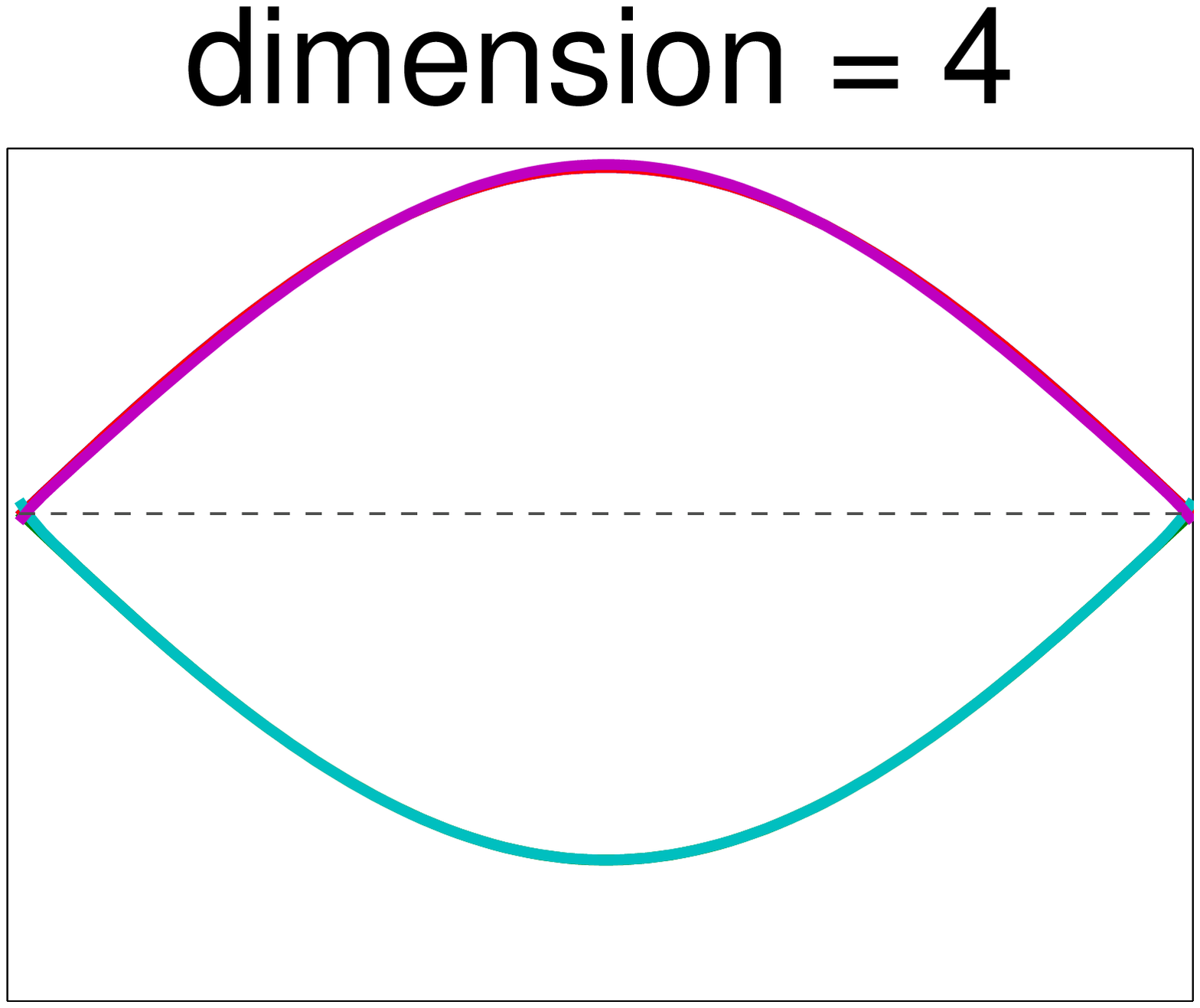} \\

\includegraphics[scale=0.12, trim = 5mm 5mm 15mm 5mm]{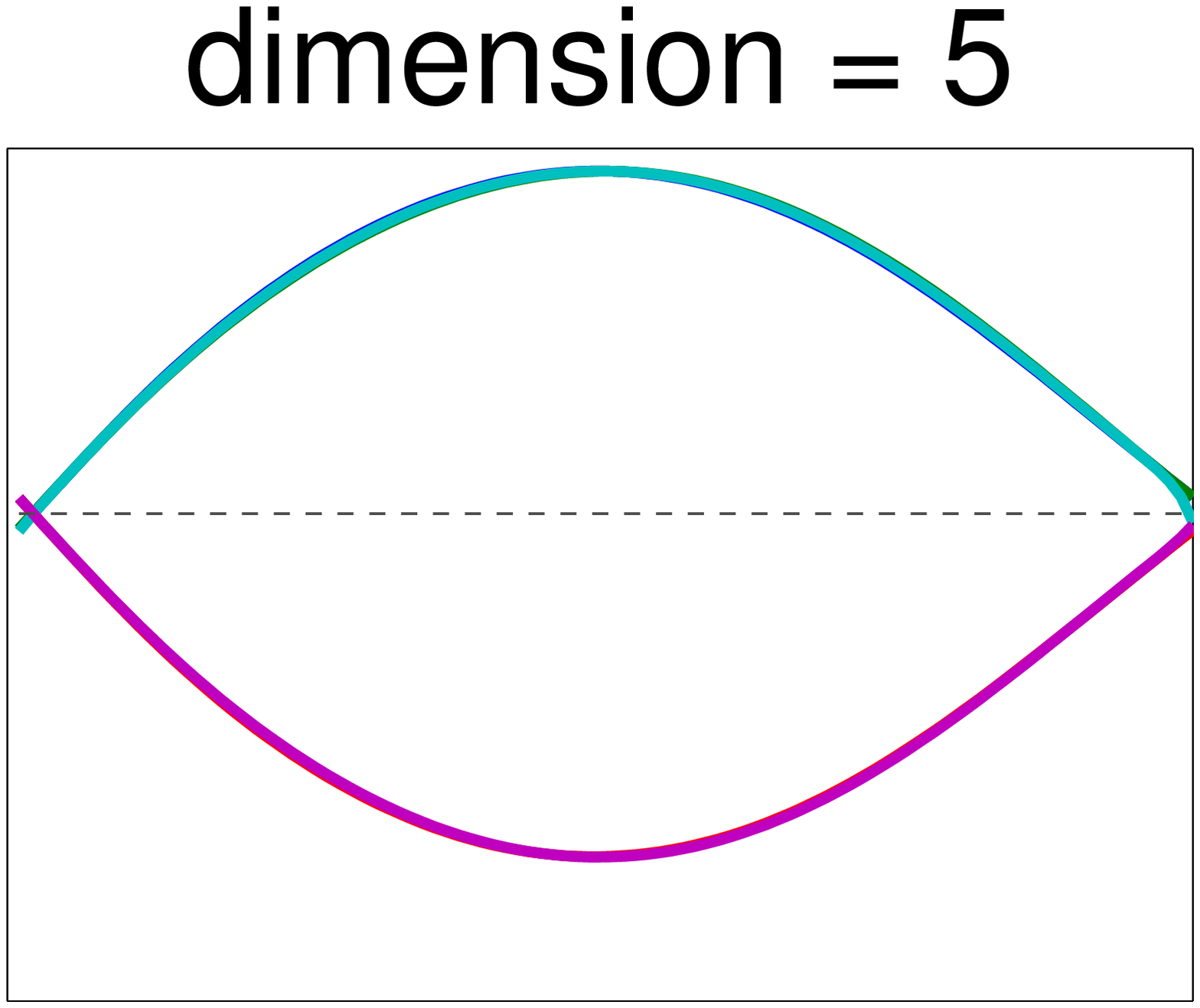} &
\includegraphics[scale=0.12, trim = 5mm 5mm 15mm 5mm]{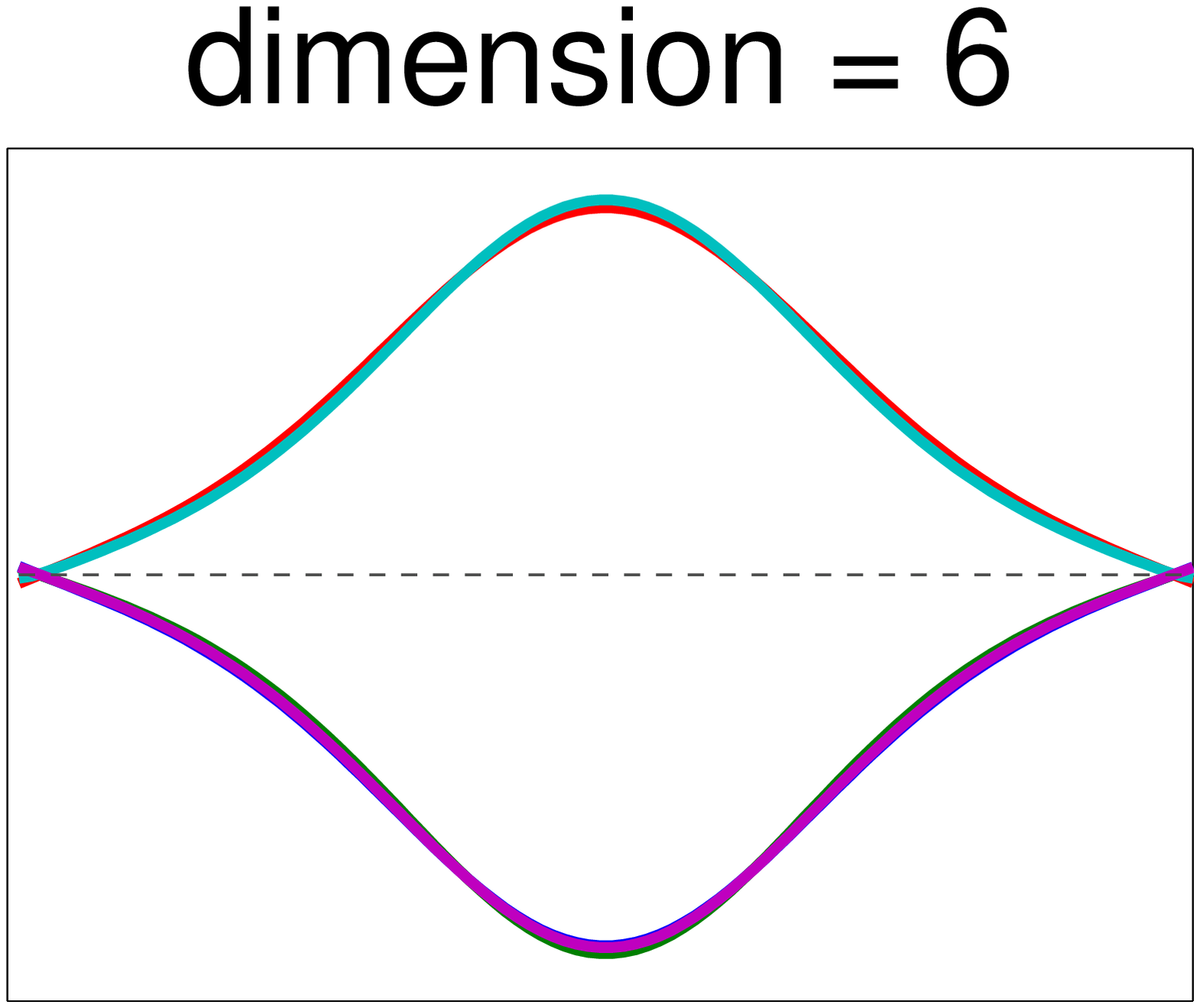} &
\includegraphics[scale=0.12, trim = 5mm 5mm 15mm 5mm]{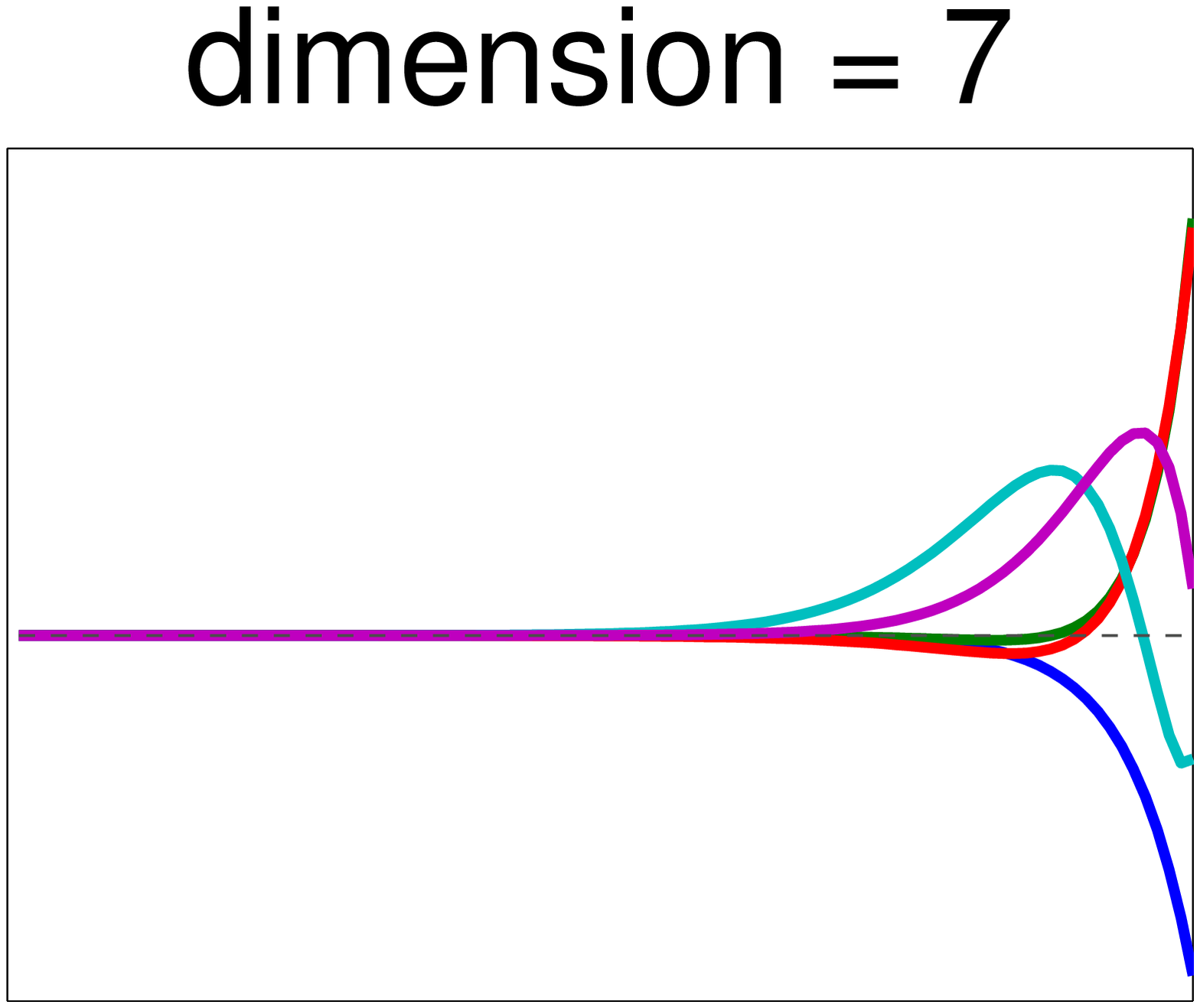} &
\includegraphics[scale=0.12, trim = 5mm 5mm 15mm 5mm]{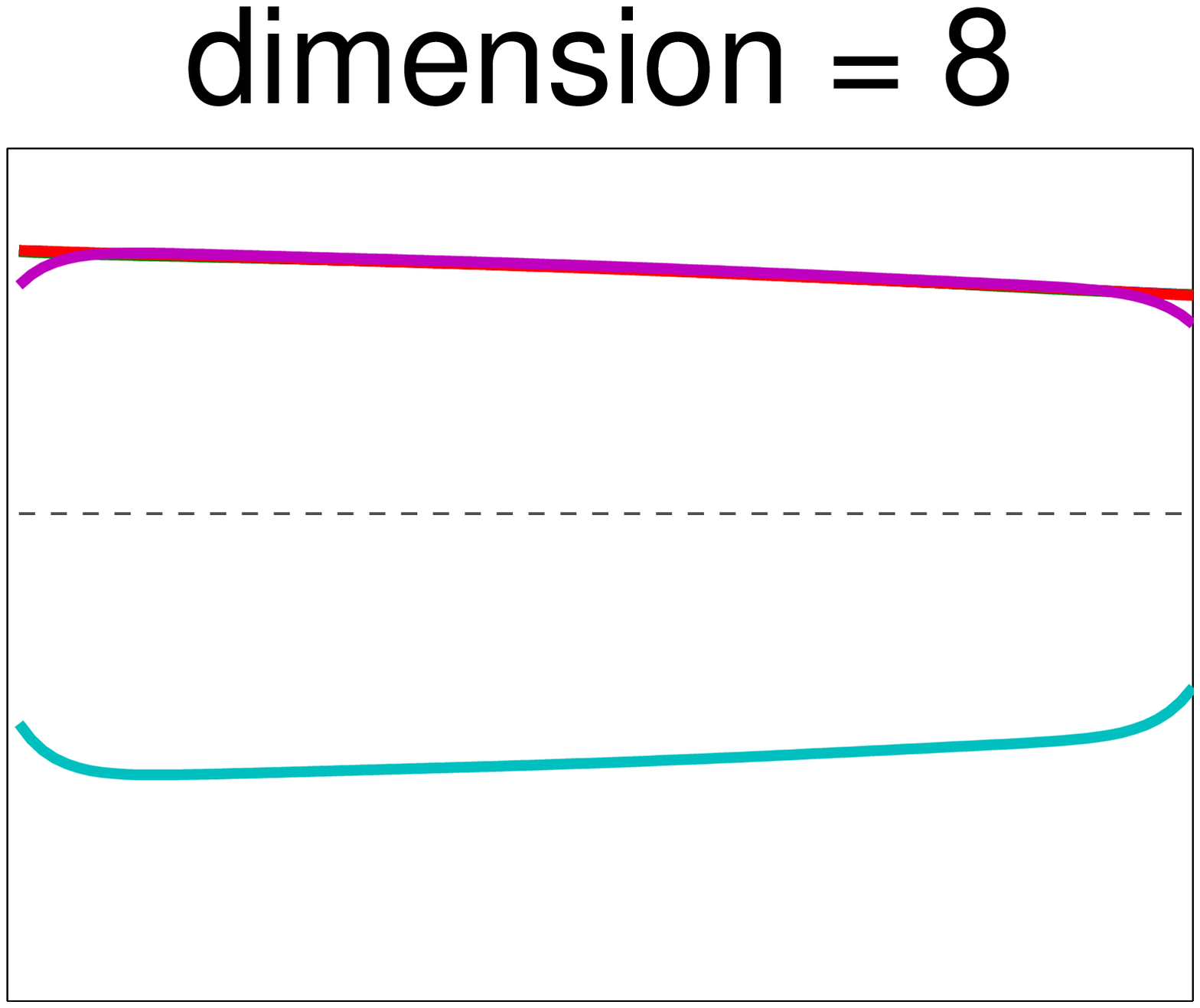} \\

\includegraphics[scale=0.12, trim = 5mm 5mm 15mm 5mm]{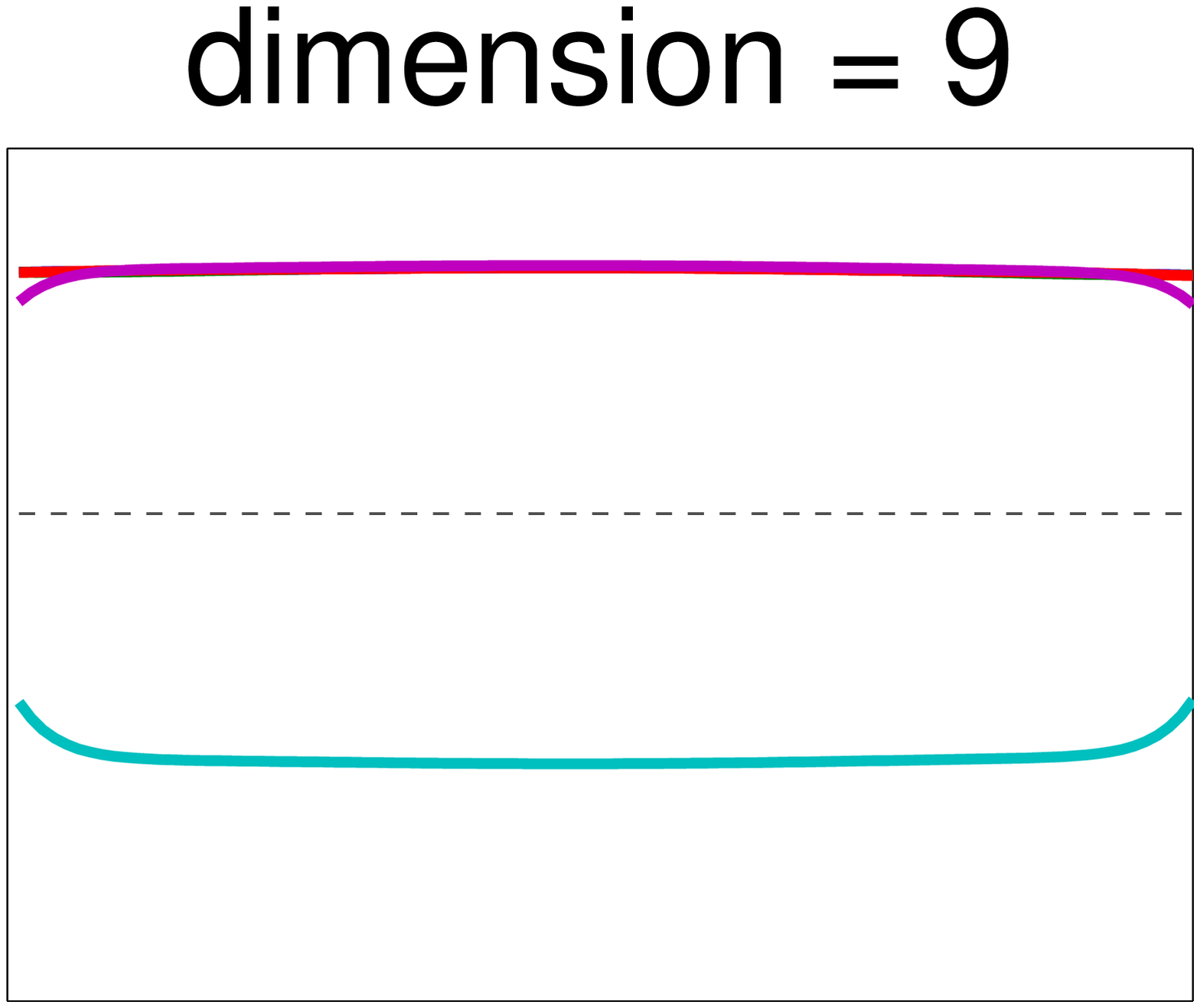} &
\includegraphics[scale=0.12, trim = 5mm 5mm 15mm 5mm]{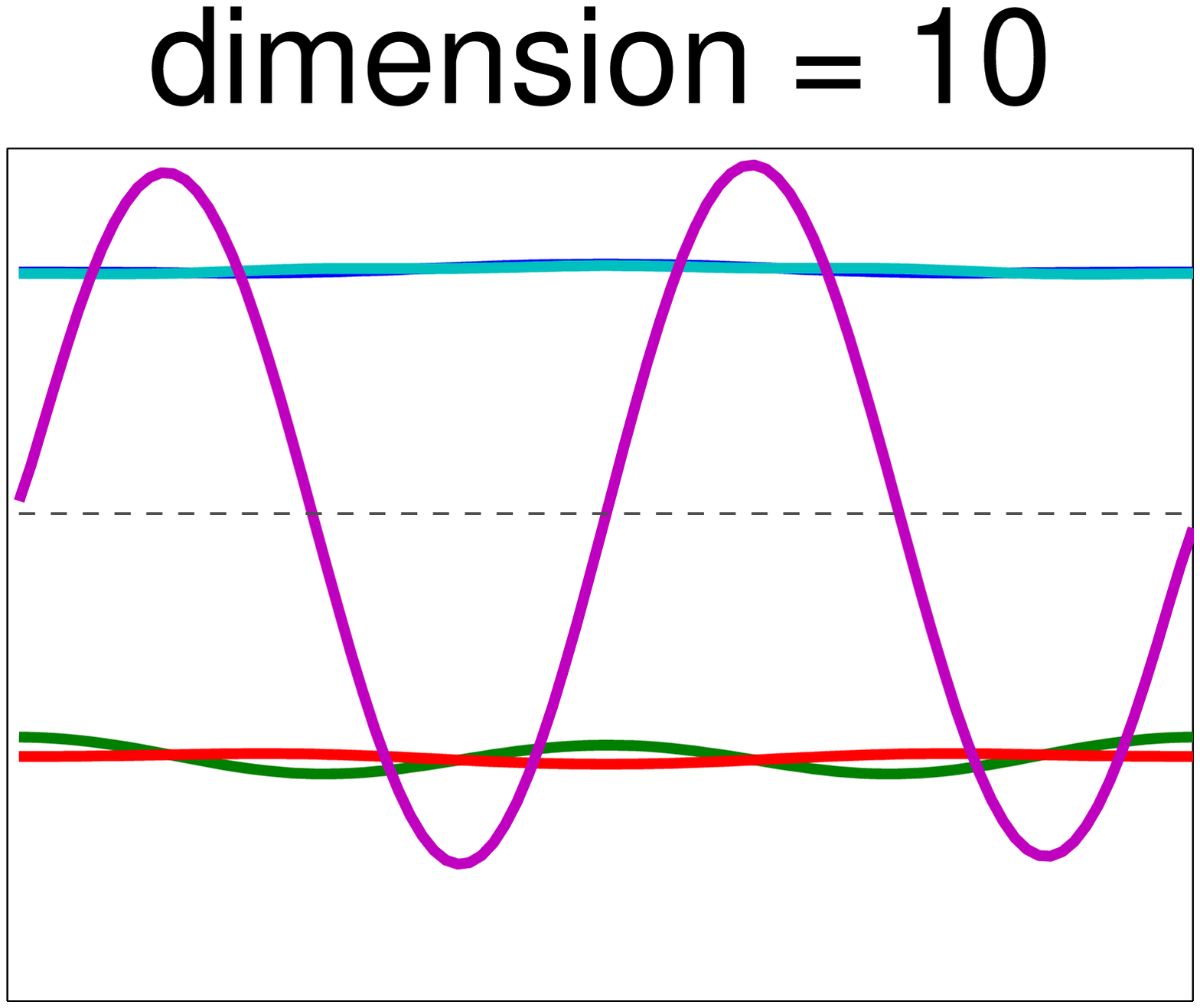} &
\includegraphics[scale=0.12, trim = 5mm 5mm 15mm 5mm]{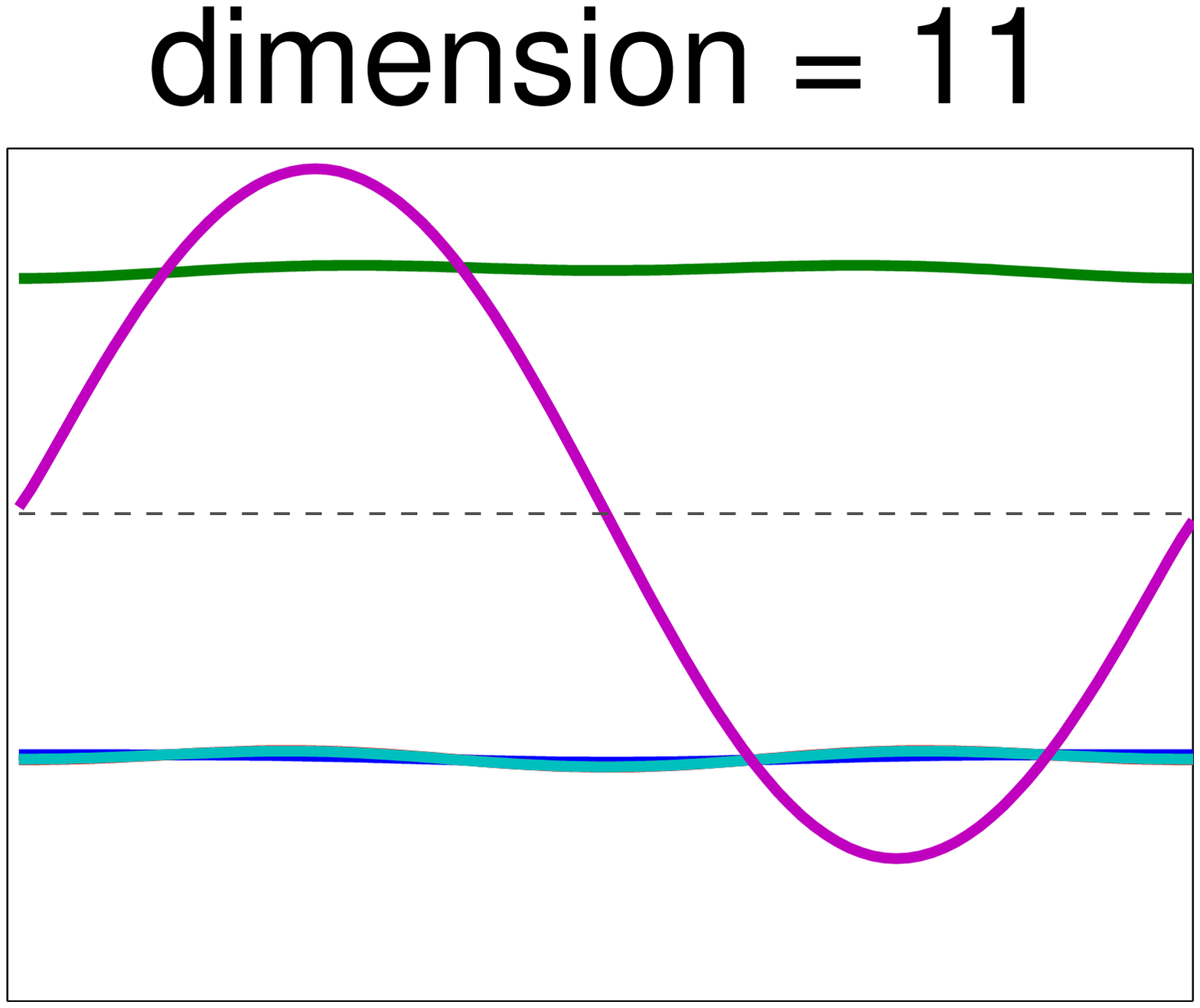} &
\includegraphics[scale=0.12, trim = 5mm 5mm 15mm 5mm]{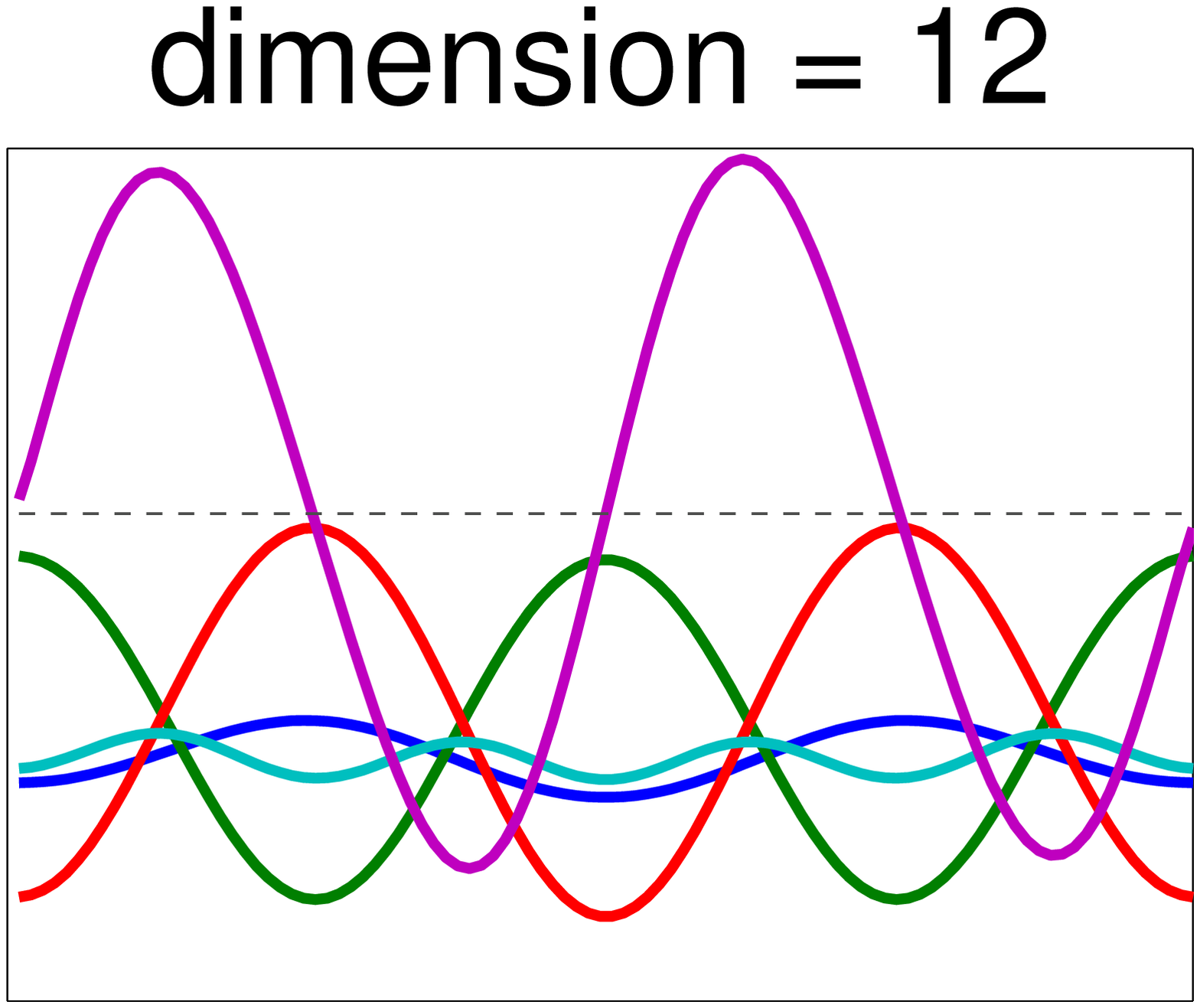}
\end{tabular}
\end{centering}
\caption{Complete basis function set for Quadcopter policy.
\label{fig:quad_basis} }
\end{figure}

\begin{figure}
\begin{centering}
\includegraphics[scale=0.23, trim = 10mm 10mm 10mm 10mm]{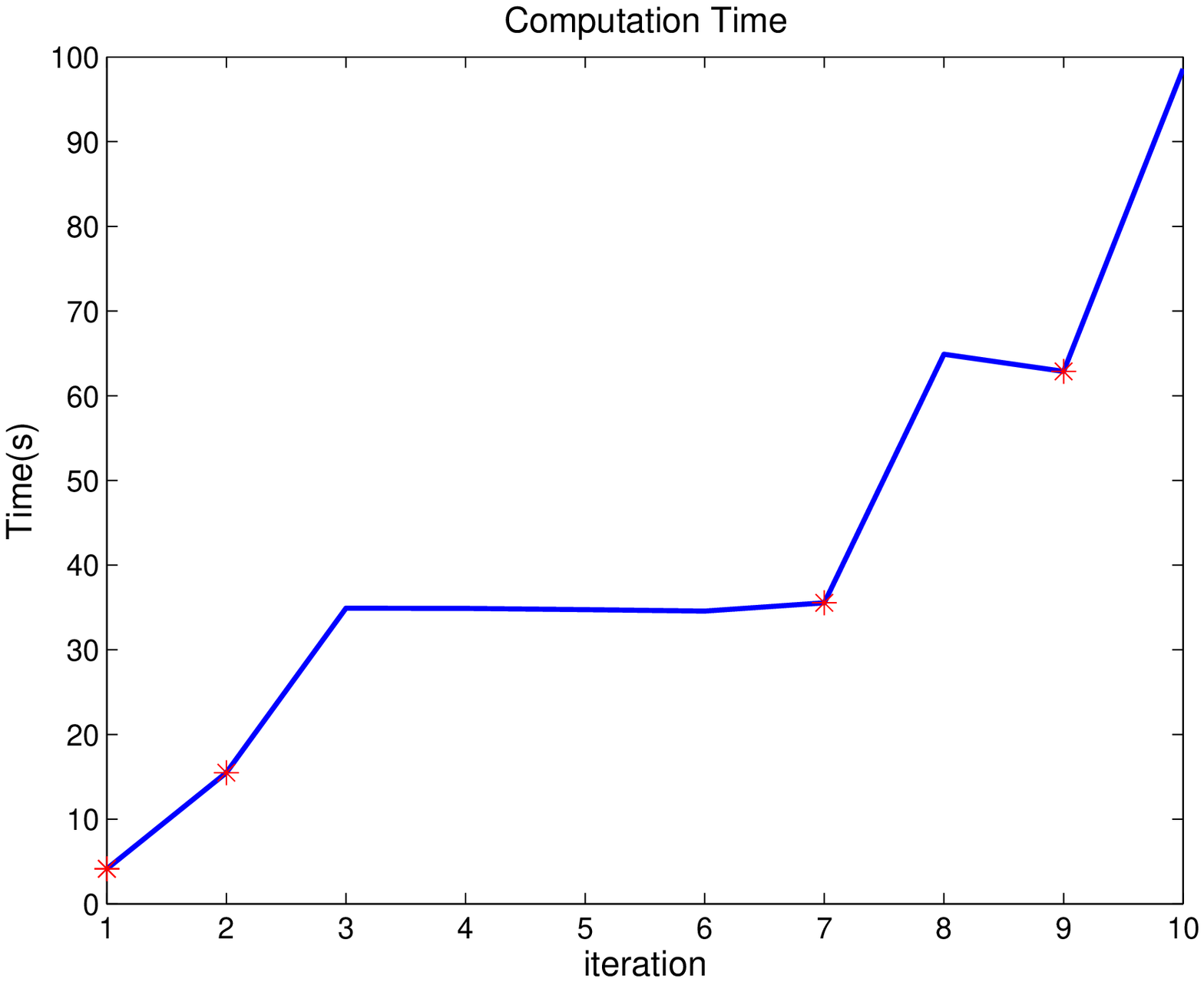}
\includegraphics[scale=0.23, trim = 10mm 10mm 10mm 10mm]{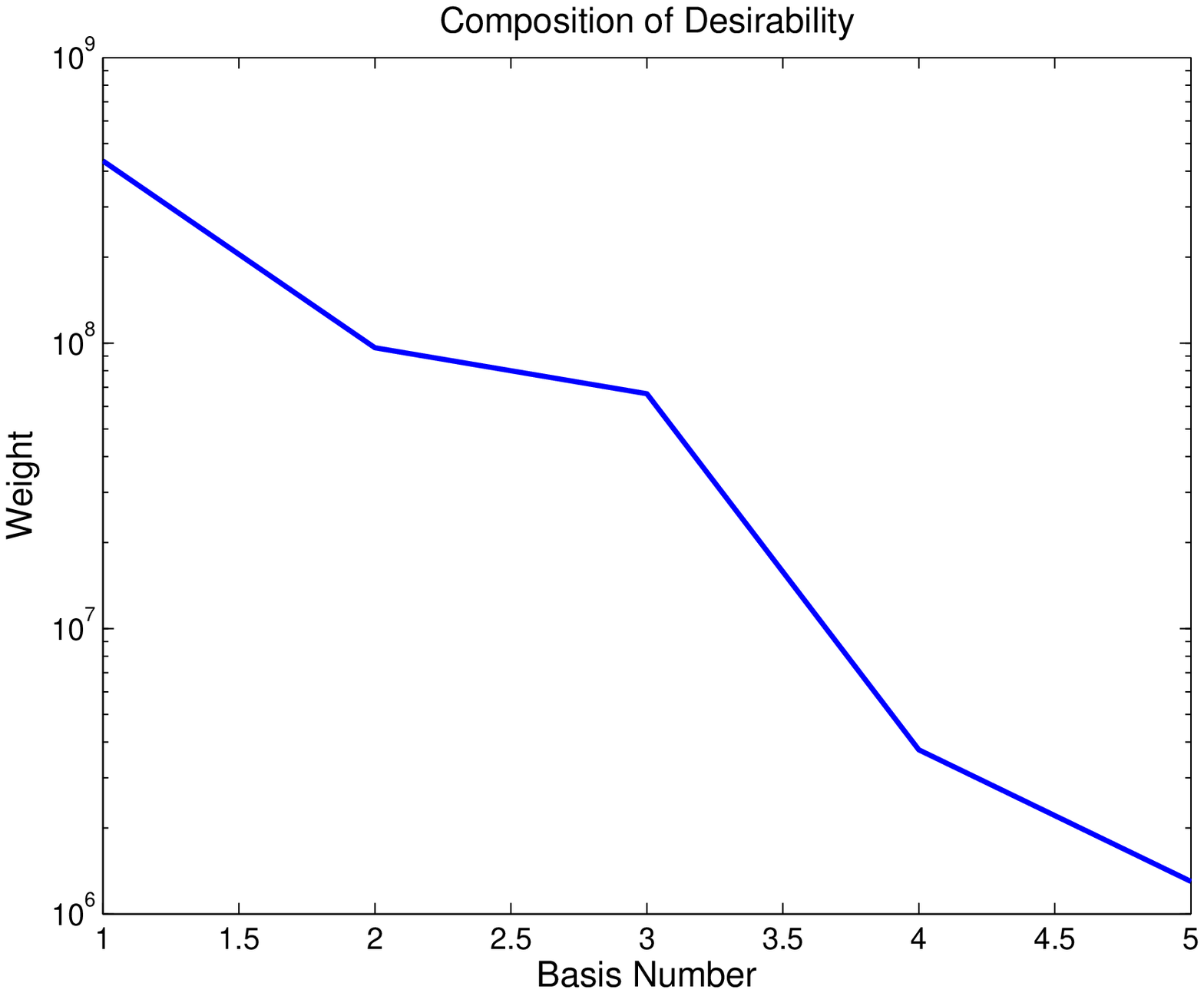}
\end{centering}
\caption{Convergence and weighting for the quadcopter solution.
\label{fig:quad_time}}
\end{figure}

\begin{figure}
\begin{centering}
\includegraphics[scale=0.5, trim = 10mm 25mm 10mm 10mm]{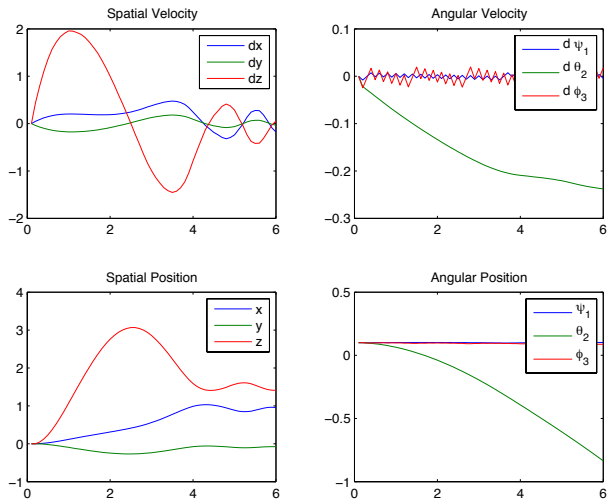}
\end{centering}
\caption{Simulation of the closed loop quadcopter system.
\label{fig:quad_sim}}
\end{figure}

Only five basis functions were computed, with the results shown in Figure \ref{fig:quad_basis}. The time for each ALS iteration is shown in Figure \ref{fig:quad_time}, along with the weighting upon each basis function. The total computation time was approximately ten minutes. Finally, Figure \ref{fig:quad_sim} shows a trajectory of the closed loop system.

\section{Discussion}


There are a number of immediate implications of this work. The first
is in the control of nonlinear distributed systems. In these problems,
additional systems manifest as additional dimensions for the PDE.
Formally, the complexity therefore grows linearly with the number of
sub-systems. As well, if the coupling between such subsystems is sparse,
it is expected that this interconnection could be simply described,
leading to low separation rank necessary to describe the coupled dynamics. 


The techniques that have been developed which rely on Sums of Squares programming \cite{Horowitz:2014tu} have been limited in degree and dimensionality due to the factorial growth in monomial basis. However,
returning to the development of the separated representation, each
rank-1 term corresponds to a single monomial. By limiting the basis
to those with high representative power, such problems may be scaled
to arbitrarily high degree and dimensionality. 

A key limitation of this work is that it requires the structural assumptions
of (\ref{eq:noise-assumption}) to obtain a linear set of equations
for which ALS may be applied. The general nonlinear value function
may not be directly solved. However, it has been shown that iterative
linearization of the nonlinear equations may be constructed in such a manner as to solve the more general HJB problem without our structural assumptions \cite{leake1967construction}. 

As alluded to in the introduction, these linear PDEs have a discrete counterpart in linearly solvable MDPs \cite{Todorov:2009je,Todorov:2006tq}. In general, MDPs must be solved through an iterative maximization process known as value or policy iteration. However, by assuming a similar restriction on the noise
of the system, specifically that it enters into the system along the same transitions actuated by the control input, Todorov has demonstrated that average cost, first exit, and finite horizon optimal control problems may be solved through a set of linear equations. It remains to be seen if the separated representation approach may also be adapted for linear MDPs. 

\subsection{Applications of the Hamilton Jacobi Bellman Solution}

The Hamilton Jacobi Bellman equation yields the optimal solution to a general form of control problem, and its impact is present in many components of control theory. Of course, the most straightforward application is that emphasized in the previous development, that of trajectory generation. The most likely trajectory of the system is in fact related to the desirability, and can be calculated from the HJB solution \cite{Todorov:2011tm}. Furthermore, although the HJB solution provides optimal trajectories, by (\ref{eq:optimal-u}) the method also provides an optimal feedback controller. The result is an architecture that is both robust and far-sighted, with the feedback controller and planner both accounting for the other. This controller has several appealing properties. In contrast to MPC-based schemes, no online computation is required, and can be seen as the optimal, continuous limit of gain scheduling. 

The ability to solve these problems for arbitrary dimension, this opens a new synthesis technique for a number of difficult problems. The first of these is the generation of Control Lyapunov Functions, which may be done by placing an exit with zero cost at the origin for the first-exit problem. The benefits of such automatic generation techniques may be seen in works such as \cite{Ames:2014jm}, where significant effort goes towards generating CLFs for particular applications, and further effort is used towards bringing these CLFs towards optimality.

\section*{Conclusion}

In this work a method to solve the Hamilton Jacobi Bellman equation
for nonlinear, stochastic systems with complexity the scales linearly
with dimension has been proposed. Although several structural assumptions
are required, systems that do not meet these may be approximated by
the introduction of noise and control effort with arbitrary magnitudes.
The implications are vast, as the curse of dimensionality no longer
necessarily prevents the use of optimal control on complex, realistic
systems. As the Hamilton Jacobi Bellman equations touch every aspect
of control theory, the techniques here hold promise in a wide variety of topics.
In particular, there are a number of important linear PDEs in control theory and estimation, including the Fokker Planck, Duncan-Mortensen-Zakai, and other equations. With the methods presented here, recourse to linearization techniques for these problems is no longer the only possibility.

\bibliographystyle{abbrv}
\bibliography{Linear_HJB,Sparse_Tensor,Moments}

\end{document}